\documentclass[11pt]{article} 
\usepackage{amsmath}
\usepackage{amssymb}
\usepackage{amsthm} 
\usepackage{mathrsfs}  

\usepackage[margin=2.0cm]{geometry} 

\usepackage[pagebackref=true]{hyperref}
\usepackage{color} 
\usepackage{enumerate}

\newtheorem{lemma}{Lemma}[section]%
\newtheorem{theorem}[lemma]{Theorem}%
\newtheorem{proposition}[lemma]{Proposition}%
\newtheorem{hypothesis}[lemma]{Hypothesis}%
\newtheorem{remark}[lemma]{Remark}%
\numberwithin{equation}{section}

\begin{document}
\title{Large maximal subgroups of almost simple classical groups}

\author{
\parbox{\textwidth}{
\centering
Fu-Gang Yin$^{\rm a}$, Ting Lan$^{\rm b,}$\footnotemark, Weijun Liu$^{\rm b,c}$, Oujie Chen$^{\rm a}$\\[2pt]
\small
\hangindent=2em\hangafter=1
$^{\rm a}$ School of Mathematics and Statistics, Beijing Key Laboratory of Biological Big Data and Topological Statistics, Beijing Jiaotong University, Beijing, P.R. China\\
$^{\rm b}$School of Mathematics and Statistics, Central South University, Changsha,  P.R. China\\
$^{\rm c}$College of General Education, Guangdong University of Science and Technology,  Dongguan, P.R. China}  
}

\renewcommand{\thefootnote}{\fnsymbol{footnote}}  
\footnotetext[1]{Corresponding author.  \\ 
E-mails: fgyin@bjtu.edu.cn(Fu-Gang Yin), lanting0603@163.com (Ting Lan), wjliu6210@126.com(Weijun Liu), 21271146@bjtu.edu.cn(Oujie Chen) \\
This work was supported by the Talent Fund 
of Beijing Jiaotong University (2025JBRC012), and the National Natural Science Foundation of China (12471022, 12331013, 12301461).}

\date{} 
\maketitle 

\begin{abstract} 
A proper subgroup \(H\) of a finite group \(G\) is called  large  if \(|H|^3 \ge |G|\). 
This paper characterises all the large maximal subgroups of almost simple groups whose socle is a finite simple classical group. 
For such a socle \(G_0\) and a core-free subgroup \(H_0\), define \(\mathcal{Q}(G_0,H_0)\) as the set of almost simple groups \(G\) with socle \(G_0\) that have a large maximal subgroup \(H\) satisfying \(H_0 = H \cap G_0\). 
We establish necessary and sufficient conditions on the pair $(G_0,H_0)$ for the set \(\mathcal{Q}(G_0,H_0)\) to be nonempty.

\end{abstract}

\section{Introduction}
A finite group $G$ is said to be \emph{almost simple} if its socle is a nonabelian simple group.
In other words, $G_0\leq G\leq \mathrm{Aut}(G_0)$ for some nonabelian simple group $G_0$. 
The theory of maximal subgroups of almost simple groups has made significant contributions to various fields.  
Regarding the orders of maximal subgroups of almost simple groups, what interesting conclusions can be drawn? 
In 1985, Liebeck~\cite{L1985} proved that, if $G_0$ is a classical group of dimension $n$  defined over a field of size $q$ then for a maximal subgroup $H$ of $G$, either $|H|<q^{3n}$, or $ H \cap G_0 $ has a well-described action on the natural projective module of $G_0$.
A similar result for exceptional groups was given by Liebeck and Saxl~\cite{LS1986}.

In 2015, Alavi and Burness~\cite{AB2015} introduced the concept of large subgroups: a proper subgroup  $H$  of a finite group \( G \) is said to be \emph{large} if \( |H|^3 \geq |G| \). 
In the same paper, they classified all the large maximal subgroups of simple groups, so it is natural to seek a classification of the large maximal subgroups of almost simple groups. 
Indeed, this has been done in~\cite{AB2015} for almost simple groups with alternating or sporadic socle.

Let $G$ be an almost simple group whose socle $G_0$ is a simple group of Lie type, and let $H$ be a core-free maximal subgroup of $G$, and set $H_0=H\cap G_0$.
Alavi and Burness~\cite{AB2015} determined all candidate pairs $(G_0,H_0)$ satisfying $|H_0|^3\geq |G_0|$.  
However, producing a complete list of pairs $(G,H)$ satsifying $|H|^3\geq |G|$ is intricate when $\mathrm{Out}(G_0)$ is large and complex. 
Instead, we pursue the more concise goal of determining all candidate pairs $(G_0,H_0)$ for which such an $H$ can exist. 
This approach has been successfully implemented for exceptional Lie types by Alavi, Bayat, and Daneshkhah \cite[Theorem 1.2]{AB2022}. 
In this work, we complete the picture by providing an analogous characterization for classical groups.

 
The study of large maximal subgroups of almost simple groups provides a relatively concise yet useful list of maximal subgroups, and it has been applied in several areas.  
For example, in finite geometry, it was proved in \cite[Corollary 4.7]{BLS2018} that if an almost simple group $G$ of Lie type acts transitively on the antiflags of a finite generalized quadrangle with $s \leq t$, then the point stabiliser in the socle of $G$ is large. A similar result was proved for locally $(G,2)$-transitive generalized quadrangles in \cite[Theorem 7.8]{BLS2021}. 
The list for large maximal subgroups of simple groups has also been used to analyse  almost simple primitive groups of bounded degree~\cite{YFX2023}.

A direct link to design theory is provided by a result of Zhou and Wang \cite[Lemma~6]{ZW2015}. They proved that for a \( G \)-flag-transitive 2-\((v, k, \lambda)\) design, the point stabiliser \(G_\alpha\)  satisfies \( |G_\alpha|^3 \geq \lambda |G| \), and hence is large in \( G \).   
A similar result was proved for  block-transitive $3$-$(v,k,1)$ designs~\cite{LLY2024}.
The authors are currently working on a project to classify the block-transitive 3-\((v, k, 1)\) designs admitting almost simple automorphism groups, which further motivates the present study.
 
Before presenting our results, we introduce some notation. 

Let $G$ be a finite almost simple classical group. 
In~\cite{Aschbacher}, Aschbacher defined $9$ subgroup collections $\mathcal{C}_1, \ldots, \mathcal{C}_8, \mathcal{S}$ of $G$, where $\mathcal{C}_1$-$\mathcal{C}_8$ are called \textit{geometric} subgroups, while $\mathcal{S}$ comprises non-geometric almost simple subgroups arising from some absolutely irreducible quasisimple groups.
We adopt the precise definition of the $\mathcal{C}_i$ collections used by Kleidman and Liebeck~\cite{K-Lie}.
Aschbacher~\cite{Aschbacher} proved that any maximal subgroup of $G$ is contained in a member of one of the $9$ collections, or in a small additional subgroup collection, denoted by $\mathcal{N}$ (following the notation of Burness and Giudici~\cite[Section 2.6]{BG2016}).
For a concise description of the subgroups in $\mathcal{N}$, see ~\cite[Table~5.9.1]{BG2016}.
Each collection $\mathcal{C}_i$ is divided into several subcollections, called \textit{types}.
We refer the reader to~\cite[Section 3.5]{K-Lie} and~\cite[Section 2.2]{BHRD2013} for further information. 
Our notation for the classical groups and for the types of subgroups in Aschbacher collections is consistent with that of~\cite{K-Lie}. 
Moreover, the notation $\mathrm{O}_n^{\circ}(q)$   is used only when both \(n\) and $q$ are odd (the sign $\circ$ is occasionally omitted). 


For a finite simple classical group $G_0$, 
we let \(\mathcal{H}(G_0)\) be the set of core-free subgroups \(H_0\) of \(G_0\) for which there exists an almost simple group $G$ with socle $G_0$ such that $G$ has a maximal subgroup $H$ satisfying $H_0=H\cap G_0$. 
For a pair $(G_0,H_0)$ with $H_0 \in \mathcal{H}(G_0)$, we define
\[
 \mathcal{Q}(G_0,H_0)=\{G_0\leq G \leq \mathrm{Aut}(G_0)\mid G \text{ has a  large maximal subgroup $H$ satisfying } H_0=H\cap G_0 \}, 
\] 
and set
\[ H_1 = N_{\mathrm{Aut}(G_0)}(H_0) \text{ and } G_1=G_0H_1.
\]
A key observation (Lemma \ref{lm:H1G1}) is that $H_1$ is maximal in $G_1$, and the set \(\mathcal{Q}(G_0, H_0)\) is nonempty if and only if \(H_1\) is large in \(G_1\). 
By investigating  investigating the largeness of \(H_1\)  in \(G_1\), we  obtain necessary and sufficient conditions  on the pair $(G_0,H_0)$ for the set \(\mathcal{Q}(G_0,H_0)\) to be nonempty; these conditions are summarised in the following four theorems.

 

\begin{theorem}\label{th:psl}
Let $G_0= \mathrm{PSL}_n(q)$ and $H_0 \in \mathcal{H}(G_0)$. Then $H_1$ is large in $G_1$  if and only if one of the following holds:
\begin{enumerate} [\rm (a)]
\item  $H_1 \in\mathcal{C}_1 \cup \mathcal{C}_8$.
\item  $H_1 \in\mathcal{C}_2$ and is of type $\mathrm{GL}_{n/t}(q)\wr \mathrm{S}_t$, where $t=2$, or $t=3$ and 
\[
q\in \{3, 4, 5, 7, 8, 9, 11, 13, 16, 17, 19, 23, 25, 27, 32, 49, 64, 81, 128\}.
\]  
\item  $H_1\in\mathcal{C}_3$ and is of type $\mathrm{GL}_{n/r}(q^r)$, where $r=2$, or $r=3$ and $q \in \{2,3,4,5,7, 8, 9, 11, 16, 27, 32\}$, or $r=5$ and $(n,q)=(5,2)$. 
\item  $H_1\in \mathcal{C}_5$ and is of type $\mathrm{GL}_n(q^{1/r})$, where $r=2$ or $3$. 
\item  $H_1\in \mathcal{C}_6$ and  $(G_0,H_0)$ is one of the following: 
\[
\begin{split}
& (\mathrm{PSL}_{4}(5),2^{4}.\mathrm{A}_{6} ), (\mathrm{PSL}_{3}(7),3^{2}.\mathrm{Q}_{8} ),\\ 
& (\mathrm{PSL}_{2}(q),\mathrm{S}_{4}) \text{ with $q\in \{7,17,23,31,41,47\}$},\\
 &(\mathrm{PSL}_{2}(q),\mathrm{A}_{4}) \text{ with $q\in \{5,11,13,19,37\}$}.
 \end{split}
\]   
\item $H_1\in \mathcal{S}$, and $(G_0,H_0)$ lies in Table~\ref{tab:S}.
\end{enumerate} 
\end{theorem}


\begin{theorem}\label{th:psu} 
Let $G_0= \mathrm{PSU}_n(q)$, where $n\geq 3$, and $H_0 \in \mathcal{H}(G_0)$. Then $H_1$ is large in $G_1$  if and only if one of the following holds:
\begin{enumerate} [\rm (a)]
\item  $H_1 \in\mathcal{C}_1$.
\item  $H_1 \in\mathcal{C}_2$ and is of type $\mathrm{GL}_{n/2}(q^2).2$, or of type  $\mathrm{GU}_{n/t}(q)\wr \mathrm{S}_t$ and one of the following holds: 
\begin{enumerate}[\rm (b.1)]
\item $t=2$;
\item $t=3$ and $q \in \{2, 3, 4, 5, 7, 8, 9, 11, 13, 16, 17, 19, 23, 25, 27, 29, 32, 49, 64, 81, 128\}$;
\item $4\leq t\leq 12$, and $(q,m,t)$ is as follows:
\[
 \begin{split}
&   ( 2, 1, 4 ),
  ( 2, 1, 5 ),
  ( 2, 1, 6 ),
  ( 2, 1, 7 ),
  ( 2, 1, 8 ),
  ( 2, 1, 9 ),
  ( 2, 1, 10 ),
  ( 2, 1, 11 ), \\&
  ( 3, 1, 4 ), 
  ( 3, 1, 5 ),
  ( 3, 1, 6 ),
  ( 4, 1, 4 ),
  ( 4, 1, 5 ),
  ( 5, 1, 4 ),
  ( 7, 1, 4 ),
  ( 8, 1, 4 ),
  ( 9, 1, 4 ).
\end{split}
\]
\end{enumerate}

\item  $H_1 \in\mathcal{C}_3$ and is of type $\mathrm{GU}_{n/r}(q^r)$ with $r=3$ and $q \in \{2, 3, 4, 5, 7, 8, 9, 16, 27, 32\}$. 
\item  $H_1\in \mathcal{C}_5$ and is of type $\mathrm{Sp}_n(q)$ or $\mathrm{O}_{n}^{\epsilon}(q)$ with $\epsilon\in\{+,-,\circ\}$ and $q$ odd, or $\mathrm{GU}_n(q^{1/r})$, where $r=2$ or $3$.   
\item  $H_1\in \mathcal{C}_6$ and  $(G_0,H_0)$ is one of the following:
\[
 (\mathrm{PSU}_{3}(5),3^{2}.\mathrm{Q}_{8} ),  
 (\mathrm{PSU}_{4}(3),2^{4}.\mathrm{A}_{6} ),   (\mathrm{PSU}_{4}(7),2^{4}.\mathrm{S}_{6} ). 
\] 
\item $H_1\in \mathcal{S}$, and $(G_0,H_0)$ lies in Table~\ref{tab:S}.
\end{enumerate} 
\end{theorem}


\begin{theorem}\label{th:psp}
Let $G_0= \mathrm{PSp}_n(q)$, where  $n\geq 4$ and $(n,q)\neq (4,2)$, and $H_0 \in \mathcal{H}(G_0)$. Then $H_1$ is large in $G_1$  if and only if one of the following holds:
\begin{enumerate} [\rm (a)]
\item  $H_1 \in \mathcal{C}_1 \cup \mathcal{C}_8$.
\item  $H_1 \in\mathcal{C}_2$ and  is of type $\mathrm{GL}_{n/2}(q)$, or of type $\mathrm{Sp}_{n/t}(q)\wr \mathrm{S}_t$, where $t \in \{2,3\}$, or $(n,t)=(8,4)$, or $(n,t)= (10,5)$ and $q\in \{3,4,5\}$.
  \item  $H_1 \in\mathcal{C}_3$ and is of type $\mathrm{GU}_{n/2}(q)$, $\mathrm{Sp}_{n/2}(q^2)$ or $\mathrm{Sp}_{n/3}(q^3)$. 
  \item  $H_1\in \mathcal{C}_5$ and is of type $\mathrm{Sp}_n(q^{1/r})$, where $r=2$ or $3$.
  \item  $H_1\in \mathcal{C}_6$ and  $(G_0,H_0)$ is one of the following
\[
 (\mathrm{PSp}_{4}(3),2^{4}.\mathrm{\Omega}^{-}_{4}(2) ),  
 (\mathrm{PSp}_{4}(5),2^{4}.\mathrm{\Omega}^{-}_{4}(2) ),   (\mathrm{PSU}_{4}(7),2^{4}.\mathrm{O}^{-}_{4}(2)),    (\mathrm{PSp}_{8}(3),2^{6}.\mathrm{\Omega}^{-}_{6}(2)). 
\] 
  \item $H_1\in \mathcal{S}$, and $(G_0,H_0)$ lies in Table~\ref{tab:S}.
  \item $H_1\in \mathcal{N}$, $G_0=\mathrm{PSp}_4(q)$ with $q$ even, and either $H_0=[q^4]{:}(q-1)^2$, or $(G_0,H_0)$ is one of the following
\[
  (\mathrm{PSp}_4(8), 7^2{:}\mathrm{D}_8), (\mathrm{PSp}_4(4), 5^2{:}\mathrm{D}_8), (\mathrm{PSp}_4(8), 9^2{:}\mathrm{D}_8), (\mathrm{PSp}_4(4), 17{:}4). 
\]
  \end{enumerate} 
\end{theorem}

\begin{theorem}\label{th:pso}
Let $G_0=\mathrm{P\Omega}^{\epsilon}_n(q)$, where $n\geq 7$ and $\epsilon \in \{\circ,+,-\}$, and $H_0 \in \mathcal{H}(G_0)$. Then $H_1$ is large in $G_1$  if and only if one of the following holds: 
\begin{enumerate} [\rm (a)]
\item  $H_1 \in \mathcal{C}_1 $.
\item  $H_1 \in\mathcal{C}_2$ and one of the following holds: 

\begin{enumerate}[\rm (b.1)]
	\item $H_1$ is of type  $\mathrm{GL}_{n / 2}(q). 2$ or $\mathrm{O}_{n / 2}^{\circ}(q)^{2}$. 
	
	\item $H_1$ is of type $\mathrm{O}_{n/2}^{\epsilon_{1}}(q) \wr \mathrm{S}_{2}$, where $\epsilon_{1}\in \{+,-,\circ\}$, and $\epsilon=\epsilon_1^2$ if $n/2 $ is even, while $\epsilon=(-1)^{(q-1)n/4}$ if $n/2 $ is odd.
	\item  The pair $(G_0,H_0)$ is one of the following: 
	\begin{align*}
& (\mathrm{P\Omega}_8^{+}(2),3^4.2^3.\mathrm{S}_4),\ (\mathrm{P\Omega}_8^{+}(3),[2^9].\mathrm{S}_4),\ (\mathrm{P\Omega}_{10}^{-}(2),3^5.2^4.\mathrm{S}_5),\ (\mathrm{P\Omega}_{12}^{-}(2), \mathrm{A}_5^3.2^2.\mathrm{S}_3),\\ 
&(\mathrm{P\Omega}_7(5),2^6.\mathrm{A}_7),\ (\mathrm{P\Omega}_8^{+}(5),2^6.\mathrm{A}_8),\ (\mathrm{P\Omega}_8^{+}(3), 2^{6}.\mathrm{A}_8), \ (\mathrm{P\Omega}_9(3), 2^{8}.\mathrm{A}_9), \ (\mathrm{P\Omega}_{10}^{-}(3), 2^{8}.\mathrm{A}_{10}), \\ 
&  (\mathrm{P\Omega}_{11}(3), 2^{10}.\mathrm{A}_{11}),    \ (\mathrm{P\Omega}_{12}^{-}(3), 2^{10}.\mathrm{A}_{12}),    \ (\mathrm{P\Omega}_{13}(3), 2^{12}.\mathrm{A}_{13}),\ (\mathrm{P\Omega}_{14}^{-}(3),2^{12}.\mathrm{A}_{14}).   
\end{align*} 
\end{enumerate} 
\item  $H_1 \in\mathcal{C}_3$ and is of type $\mathrm{O}_{n/2}^{\epsilon}(q^{2})$, $\mathrm{O}_{n/2}^{\circ}(q^{2})$ or $\mathrm{GU}_{n/2}(q)$. 
\item $H_1 \in\mathcal{C}_4$ and is of type $\mathrm{Sp}_{2}(q) \otimes \mathrm{Sp}_{n/2}(q)$ and $(n,\epsilon)=(8,+)$ or $(12,+)$.
\item  $H_1\in \mathcal{C}_5$ and is of type $\mathrm{O}_{n}^{\epsilon'}(q^{1/2})$ with $\epsilon'^2=\epsilon$, or $\mathrm{O}_{n}^{\epsilon}(q^{1/3})$.  
\item  $H_1\in \mathcal{C}_6$ and $(G_0,H_0)=(\mathrm{P\Omega}_8^{+}(3),2^6.\Omega_{6}^+(2))$.  
\item $H_1\in \mathcal{S}$, and $(G_0,H_0)$ lies in Table~\ref{tab:S}.
\item  $H_1\in\mathcal{N} $, $G_0=\mathrm{P\Omega}_8^{+}(q)$, and $H_0$ is of type $\mathrm{GL}_1(q)\times \mathrm{GL}_3(q)$, $\mathrm{GU}_1(q)\times \mathrm{GU}_3(q)$, $ \mathrm{G}_2(q)$, $P_{1,3,4}$ (a non-maximal parabolic subgroup), or $(G_0,H_0)\in\{(\mathrm{P\Omega}_8^{+}(2), \mathrm{D}_{10}^2.2.\mathrm{S}_2),(\mathrm{P\Omega}_8^{+}(3),[2^9].\mathrm{SL}_3(2))\}$.  
\end{enumerate} 
\end{theorem}


Our proofs rely heavily on results and techniques of Alavi and Burness in \cite{AB2015}. In particular, when \(H_0\) itself is already large in \(G_0\), the corresponding pairs are automatically included in our final list.
In the lemmas presented in our proofs, we explicitly note which pairs \((G_0, H_0)\) satisfy that \(H_1\) is large in \(G_1\) while \(H_0\) itself is not large in \(G_0\).

\section{Preliminaries}
In this section, we fix notation and present some preliminary results that will be frequently used throughout the paper.

For positive integers \(a_1,\ldots,a_m\), we write \((a_1,\ldots,a_m)\) and \([a_1,\ldots,a_m]\) for their greatest common divisor and least common multiple, respectively. 
Let $n$ be a positive integer and $p$ a prime. We write $n_p$ for the largest power of $p$ dividing $n$.

For a positive integer $n $, let $\mathrm{C}_n$, $\mathrm{D}_n$ and $\mathrm{Q}_n$  denote the cyclic group, dihedral group and quaternion group of order $n$,  respectively.
Moreover, $\mathrm{A}_n$ and $\mathrm{S}_n$  denote the alternating and symmetric group of degree $n$, respectively.
For two groups $A$ and $B$, we write $A.B$ for an extension of $A$ by $B$, and $A{:}B$ for a split extension of $A$ by $B$, and $A\times B$ for the direct product of $A$ and $B$.

We begin with the order formulas of classical groups. Table~\ref{tb:orderofsimplegroups} is taken from~\cite[Table 5.1.A]{K-Lie}, where $q=p^e$ for some prime $p$. 

\begin{table}[h]
  \caption{The orders of simple classical groups and their outer automorphism groups.}\label{tb:orderofsimplegroups}
\[
\begin{array}{ccccl}
\hline
G_0 &  
   d & |\mathrm{Out}(G_0)| & |G_0| \\
\hline
\mathrm{PSL}_n(q)  & (n,q - 1) & \begin{array}{c}
2de, n\geq3\\
de, n = 2
\end{array} & \frac{1}{d}q^{n(n - 1)/2}\prod_{i = 2}^{n}(q^{i}-1)\\
\hline 
\mathrm{PSU}_n(q)  & (n,q + 1) & \begin{array}{c}
  2de, n\geq3\\
  de, n = 2
  \end{array} & \frac{1}{d}q^{n(n - 1)/2}\prod_{i = 2}^{n}(q^{i}-(-1)^{i})\\
  \hline
\mathrm{PSp}_{2m}(q)   & (2,q - 1) & \begin{array}{c}
de, m\geq3\\
2e, m = 2
\end{array} & \frac{1}{d}q^{m^{2}}\prod_{i = 1}^{m}(q^{2i}-1)\\
\hline
\begin{array}{c}
\Omega_{2m + 1}^{\circ }(q)\\
q\text{ odd}
\end{array} & 2 & 2e & \frac{1}{2}q^{m^{2}}\prod_{i = 1}^{m}(q^{2i}-1)\\
\hline
\begin{array}{c}
  \mathrm{P\Omega}_{2m}^{+}(q)\\
m\geq3
\end{array}  &(4,q^{m}-1) & \begin{array}{c}
2de, m\neq4\\
6de, m = 4
\end{array} & \frac{1}{d}q^{m(m - 1)}(q^{m}-1)\prod_{i = 1}^{m - 1}(q^{2i}-1)\\
\hline
\begin{array}{c}
  \mathrm{P\Omega}_{2m}^{-}(q)\\
m\geq2
\end{array}  & (4,q^{m}+1) & 2de & \frac{1}{d}q^{m(m - 1)}(q^{m}+1)\prod_{i = 1}^{m - 1}(q^{2i}-1)\\
\hline
\end{array}
\] 
\end{table}

The next two propositions are derived from~\cite[Table 2.1.C]{K-Lie} and~\cite[Proposition 2.9.1]{K-Lie}.

\begin{proposition}\label{prop:orders}
Let $d$ be as in Table~\ref{tb:orderofsimplegroups}.
\begin{enumerate} [\rm (a)]
  \item $|\mathrm{SL}_n(q)|=|\mathrm{PGL}_n(q)|=d|\mathrm{PSL}_n(q)|$, and $|\mathrm{GL}_n(q)|=(q-1) |\mathrm{SL}_n(q)| $.
  \item $|\mathrm{SU}_n(q)|=|\mathrm{PGU}_n(q)|=d|\mathrm{PSU}_n(q)|$, and $|\mathrm{GU}_n(q)|=(q-1) |\mathrm{SU}_n(q)| $.
  \item $|\mathrm{Sp}_{2m}(q)|  =d|\mathrm{PSp}_{2m}(q)|$.
  \item If $q$ is odd, then $|\mathrm{\Omega}_{2m+1}(q)|=|\mathrm{P\Omega}_{2m+1}(q)|$, and $|\mathrm{SO}_{2m+1}(q)|=2|\mathrm{\Omega}_{2m+1}(q)|$, and $|\mathrm{O}_{2m+1}(q)|=2|\mathrm{SO}_{2m+1}(q)|$.
  \item For every $\epsilon\in \{+,-\}$ and $m\geq 1$, $|\mathrm{PO}_{2m}^{\epsilon}(q)|=|\mathrm{SO}_{2m}^{\epsilon}(q)|=2|\mathrm{\Omega}_{2m}^{\epsilon}(q)|=(2,q-1)|\mathrm{PSO}_{2m}^{\epsilon}(q)|$, and $|\mathrm{O}_{2m}^{\epsilon}(q)|=(2,q-1)|\mathrm{SO}_{2m}^{\epsilon}(q)|$, and  
  \[
   |\mathrm{\Omega}_{2m}^{\epsilon}(q)| =\frac{1}{(2,q-1)}q^{m(m - 1)}(q^{m}+1)\prod_{i = 1}^{m - 1}(q^{2i}-1).
  \] 
\end{enumerate}
\end{proposition}

\begin{proposition}\label{prop:iso}
\begin{enumerate}[\rm (a)]
  \item $\mathrm{PSL}_2(q)\cong \mathrm{PSU}_2(q) \cong \mathrm{PSp}_2(q) \cong \mathrm{P\Omega}_3 (q) $.
  \item $\mathrm{\Omega}^{+}_2(q)\cong \mathrm{C}_{ \frac{q-1}{(2,q-1)}}$,  $\mathrm{\Omega}^{-}_2(q)\cong \mathrm{C}_{  \frac{q+1}{(2,q-1)}}$,  $\mathrm{\Omega}^{+}_4(q)\cong \mathrm{SL}_{2}(q)\circ \mathrm{SL}_{2}(q)$,  $\mathrm{\Omega}^{-}_4(q)\cong \mathrm{PSL}_{2}(q^2)$, $\mathrm{P\Omega}_5 (q)\cong \mathrm{PSp}_4(q)$, $\mathrm{P\Omega}^{+}_6(q)\cong \mathrm{PSL}_{4}(q) $,  $\mathrm{P\Omega}^{-}_6(q)\cong \mathrm{PSU}_{4}(q)$, and  $\mathrm{P\Omega}_{2m+1}(q)\cong \mathrm{PSp}_{2m}(q)$ for even $q$. 
\end{enumerate}
\end{proposition} 

We now present some bounds on the orders of classical groups, which are taken from~\cite{AB2015}.

\begin{lemma}[{\cite[Lemma~4.2]{AB2015}}] \label{lm:LUSO}
Let $q$ be a prime power.
\begin{enumerate}[\rm (i)]
\item If $n,q\geq 2$ then 
\begin{align*}
 &(1-q^{-1}-q^{-2} )q^{n^{2}} <|\mathrm{GL}_{n}(q) |\leq (1-q^{-1} ) (1-q^{-2} )q^{n^{2}},\\ 
 &(1+q^{-1} ) (1-q^{-2} )q^{n^{2}} <|\mathrm{GU}_{n}(q) |\leq (1+q^{-1} ) (1-q^{-2} ) (1+q^{-3} )q^{n^{2}}. 
\end{align*}

\item  If $n\geq4$ and $q\geq2$ then 
\[(1-q^{-2}-q^{-4})q^{\frac{1}{2}n(n+1)}<|\mathrm{Sp}_{n}(q)|\leq (1-q^{-2})( 1-q^{-4}) q^{\frac{1}{2}n(n+1)}.\]

\item  If $n\geq 5$ and $q\geq2$ then 
\begin{align*}
&(1-q^{-2}-q^{-4})q^{\frac{1}{2}n(n-1)} <|\mathrm{SO}_n^\circ(q)| <(1-q^{-2})(1-q^{-4})q^{\frac{1}{2}n(n-1)},\\ 
&(1-q^{-2}-q^{-4})(1-q^{-\frac{n}{2}})q^{\frac{1}{2}n(n-1)} <\frac{|\mathrm{SO}_n^+(q)|}{(2,q)} <(1-q^{-2})(1-q^{-4})q^{\frac{1}{2}n(n-1)},\\  
&(1-q^{-2}-q^{-4})q^{\frac{1}{2}n(n-1)} <\frac{|\mathrm{SO}_n^-(q)|}{(2,q)} <(1-q^{-2})(1-q^{-4})(1+q^{-\frac{n}{2}})q^{\frac{1}{2}n(n-1)}.  
\end{align*} 
\end{enumerate}
\end{lemma}

\begin{lemma}[{\cite[Corollary~4.3]{AB2015}}] \label{lm:LUSO1} 
\begin{enumerate}[\rm (i)]
\item If $n\geqslant2$ then
\[ 
q^{n^2-2}<| \mathrm{PSL}_n(q)|\leq (1-q^{-2})q^{n^2-1}.
\]

\item If $n\geq 3$ then
\[
(1-q^{-1})q^{n^2-2}<| \mathrm{PSU}_n(q)|\leq (1-q^{-2})(1+q^{-3})q^{n^2-1}<q^{n^2-1}.
\]

\item If $n\geq 4$ then
\[
\frac{1}{2(2,q-1)}q^{\frac{1}{2}n(n+1)}< |\mathrm{PSp}_{n}(q) |<q^{\frac{1}{2}n(n+1)}.
\]  

\item If $n\geq 7$ then
\[
\frac{1}{4(2,n)}q^{\frac{1}{2}n(n-1)}< |\mathrm{P}\Omega_{n}^{\epsilon}(q) |<q^{\frac{1}{2}n(n-1)}.
\]   
 \end{enumerate}

In particular, $|G|>\frac{1}{8}q^{\frac{1}{2}n(n-1)}$ for every finite simple classical group $G$ of dimension $n$ over $\mathbb{F}_{q}$.
\end{lemma}
By Propositions~\ref{prop:orders},~\ref{prop:iso} and Lemma~\ref{lm:LUSO}, we can prove the following result.  
  
\begin{lemma}\label{lm:Omega}
Let $\epsilon\in \{+,-,\circ \}$ and $n\geq 2$. Then $|\Omega_n^\epsilon(q)|<\frac{1}{(2,q-1)}q^{\frac{n(n-1)}{2}} $, except when $n=2$ and $\epsilon=-$.
  Moreover, $|\Omega_2^-(q)|= (q+1)/{(2,q-1)}$.
\end{lemma}

The following lemma is a consequence of the AM-GM inequality, which states that \( (a_1+ \cdots+a_n)/n\geq \sqrt[n]{a_1 \cdots a_n} \).
\begin{lemma}\label{lm:n!}
  Let \( t\geq 2 \) be a positive integer. Then \( t!<(\frac{t+1}{2})^t=2^{t\log_2(\frac{t+1}{2})} \).
\end{lemma}

The next result is straightforward to prove.
\begin{lemma}\label{lm:q2f}
Let \( p \) be a prime and let \( q=p^e \), where \( e \) is a positive integer. 
Then \( e^2=\log_p(q)^2 \leq q \), except when \( q=8 \).
In particular, \( \log_p(q)^2 \leq 9q/8\) for any \( q \). 
\end{lemma}

\section{Proof of Theorems~\ref{th:psl}--\ref{th:pso}}\label{sec:proof} 
Throughout this section,  $G_0$ is a finite simple classical group, $H_0\in \mathcal{H}(G_0)$, and $q$ is a power of a  prime $p$.

\subsection{The basic idea} \label{sec:proof-1}

Recall the set 
\begin{align*}
& \mathcal{Q}(G_0,H_0)=\{G_0\leq G \leq \mathrm{Aut}(G_0)\mid G \text{ has a  large maximal subgroup $H$ satisfying } H_0=H\cap G_0 \},\\
&H_1 = N_{\mathrm{Aut}(G_0)}(H_0) \text{ and } G_1=G_0H_1.
\end{align*} 
Our goal is to determine those pairs $(G_0,H_0)$ such that $\mathcal{Q}(G_0,H_0)$ is nonempty. 

Throughout Section~\ref{sec:proof}, we set  
\[
\mathcal{O}_{1}=G_{1}/G_0 \leq \mathrm{Out}(G_0).
\]
Note that, $G_{1}=G_0H_{1}$ implies $H_{1}/H_0\cong  G_{1}/G_0$, and so $G_1=G_0.\mathcal{O}_1$ and  $H_1=H_0.\mathcal{O}_1$. 

\begin{lemma}\label{lm:H1G1}
The group $H_{1}$ is maximal in $G_{1}$.  
Moreover, the set $\mathcal{Q}(G_0,H_0)$ is nonempty if and only if $H_{1}$ is large in $G_{1}$. 
\end{lemma}

\begin{proof}
Suppose for a contradiction that $H_{1}$ is not maximal in $G_{1}$. Then there exists a subgroup $M$ of $G_{1}$ such that $H_{1} < M < G_{1}$. 
Since $H_0\in \mathcal{H}(G_0)$, there exists an almost simple group $G$ with socle $G_0$ such that  $G$ has a maximal subgroup $H$ satisfying $H\cap G_0=H_0$. 
Since $H_0=H\cap G_0 \unlhd H$, we have $H\leq H_1<M$.
Note that $M\cap G \geq H$ and $H$ is maximal in $G$.
Therefore, $M\cap G=H$, and so $M \cap G_0 =H \cap G_0 = H_0$. 
Now, because $G_0$ is normal in $G_{1}$, we have $M \cap G_0 = H_0$ is normal in $M$. This implies $M \leq N_{G_{1}}(H_0) = H_{1}$. However, this contradicts the assumption $H_{1} < M$.  
Therefore, $H_{1}$ must be maximal in $G_{1}$. 

If $H_{1}$ is large in $G_{1}$, then $G_1\in \mathcal{Q}(G_0,H_0)$ and hence $\mathcal{Q}(G_0,H_0)$ is nonempty. 

Now, suppose that $\mathcal{Q}(G_0,H_0)$ is nonempty. 
Then there exists some $G\in \mathcal{Q}(G_0,H_0)$ such that $G$ has a large maximal subgroup $H$ satisfying $H_0=H\cap G_0$.
Write $G=G_0.\mathcal{O}$ with $\mathcal{O}\leq \mathrm{Aut}(G_0)$.  
Then $H=N_{G}(H_0)\leq H_{1}$ and $H=H_0.\mathcal{O}$.
Note that $H_{1}/H_0\cong G_{1}/G_0=\mathcal{O}_{1}$ and $\mathcal{O}\leq \mathcal{O}_{1}$ (as $H_0.\mathcal{O}=H\leq H_{1}=H_0.\mathcal{O}_1$).
Thus, 
\[
\frac{|H_{1}|^3}{|G_{1}|}=\frac{|H_0|^3}{|G_0|}\cdot  |\mathcal{O}_{1}|^2 \geq \frac{|H_0|^3}{|G_0|}\cdot  |\mathcal{O}|^2=\frac{|H_0.\mathcal{O}|^3}{|G_0.\mathcal{O}|}=\frac{|H|^3}{|G|}\geq 1.
\]
It follows that $H_{1}$ is large in $G_{1}$, and so $G_{1} \in \mathcal{Q}(G_0,H_0)$.
\end{proof}

According to Lemma~\ref{lm:H1G1},  to  prove Theorems~\ref{th:psl}--\ref{th:pso}, it is sufficient to check whether 
  $H_{1}$ is large in $G_{1}$. This is equivalent to checking  the following inequality:
\begin{equation}\label{eq:H1}
|G_0|\leq   \vert H_0\vert^3  \vert \mathcal{O}_{1}  \vert^2 .
\end{equation} 
 
We note that the group $\mathcal{O}_{1}$ is consistent with the ``stabilizer'' column in~\cite[Table 3.5.G]{K-Lie} (by~\cite[(3.2.3)]{K-Lie}), and also the ``Stab'' columns in the tables of~\cite[Section 8]{BHRD2013}. This describes the stabiliser of one $G_0$-conjugacy class of subgroups $H_0$ under the action of $\mathrm{Out}(G_0)$.
Let 
\begin{equation}\label{eq:c}
  c=|\mathrm{Out}(G_0)|/|\mathcal{O}_{1}|.  
\end{equation}
 
Then $c$ represents the number of $G_0$-conjugacy classes of subgroups $H_0$. Note that, by Aschbacher's theorem (see~\cite[Theorem 3.1.1]{K-Lie}), $\mathrm{Out}(G_0)$ acts transitively on those $c$ classes.


\subsection{A technique of Alavi and Burness}\label{sec:techniques}
In this subsection,  we recall a key technique in~\cite{AB2015} for handling the boundary case where \( |H_0|^3 \approx |G_0| \). This method will be used repeatedly in our subsequent analysis.


Consider, for instance, the case where \( G_0 = \mathrm{PSL}_n(q) \), and \( H_0 \) is a \(\mathcal{C}_2\)-subgroup of type \( \mathrm{GL}_m(q) \wr \mathrm{S}_3 \) with \( n = 3m \). Here, we have \( |G_0| = d^{-1}(q-1)^{-1}|\mathrm{GL}_{3m}(q)| \), \( |H_0| = d^{-1}(q-1)^{-1}|\mathrm{GL}_{m}(q)|^3 \cdot 3! \), and \( |\mathcal{O}_{1}| = 2d\log_p(q) \), where \( d = (n, q-1) \). Now, 
\[
\frac{|H_0|^3|\mathcal{O}_{1}|^2}{|G_0|} = \frac{(3!)^3 \cdot (4d^2\log_p(q)^2)}{d^2(q-1)^2} \cdot \frac{|\mathrm{GL}_{m}(q)|^9}{|\mathrm{GL}_{3m}(q)|} = \frac{864\log_p(q)^2}{(q-1)^2} \cdot \frac{|\mathrm{GL}_{m}(q)|^9}{|\mathrm{GL}_{3m}(q)|}.
\]
Then Eq.~\eqref{eq:H1} becomes 
\begin{equation} \label{eq:hq}
  \frac{|\mathrm{GL}_{m}(q)|^9}{|\mathrm{GL}_{3m}(q)|} \geq \frac{(q-1)^2}{864\log_p(q)^2} := h(q).  
\end{equation} 

By Lemma~\ref{lm:LUSO}, \( |\mathrm{GL}_{3m}(q)| \approx q^{9m^2} \approx |\mathrm{GL}_{m}(q)|^9 \), along with precise bounds:
\begin{align*}
& (1-q^{-1}-q^{-2})q^{9m^{2}} < |\mathrm{GL}_{3m}(q)| < (1-q^{-1})(1-q^{-2})q^{9m^{2}}, \\
& (1-q^{-1}-q^{-2})^9 q^{9m^{2}} < |\mathrm{GL}_{m}(q)|^9 < (1-q^{-1})^9(1-q^{-2})^9 q^{9m^{2}}.    
\end{align*}
This yields the inequality  
\[
f(q) := \frac{(1-q^{-1}-q^{-2})^9}{(1-q^{-1})(1-q^{-2})} < \frac{|\mathrm{GL}_{m}(q)|^9}{|\mathrm{GL}_{3m}(q)|} < \frac{(1-q^{-1})^9(1-q^{-2})^9}{1-q^{-1}-q^{-2}} := g(q).
\]  
Since \( (1-q^{-1})^2 = 1-2q^{-1}+q^{-2} = 1-q^{-1}-(q^{-1}-q^{-2}) \leq 1-q^{-1}-q^{-2} \), we have \( g(q) < 1 \). Then Eq.~\eqref{eq:hq} implies \( 1 > h(q) \), and hence
\[
(q-1)^2 < 864\log_p(q)^2.
\]
Direct computation reveals that the above inequality holds precisely when 
\[
q \in \{2, 3, 4, 5, 7, 8, 9, 11, 13, 16, 17, 19, 23, 25, 27, 29, 32, 49, 64, 81, 128\}.
\] 
For these values of \( q \), further computation shows that 
\( h(q) < g(q) \) holds if and only if \( q \notin \{2, 29\} \), and moreover, \( h(q) < f(q) \) always holds for \( q \notin \{2, 29\} \). 
Notice that, for a given \( q \), if Eq.~\eqref{eq:H1} holds, then we must have \( h(q) < g(q) \), and if \( h(q) < f(q) \) holds, then Eq.~\eqref{eq:H1} must hold.
Hence, \( H_{1} \) is large in \( G_{1} \) if and only if 
\[
q \in \{ 3, 4, 5, 7, 8, 9, 11, 13, 16, 17, 19, 23, 25, 27, 32, 49, 64, 81, 128\}.
\] 



\subsection{Geometric subgroups of linear groups} \label{sec:proof-psl}
In this subsection, we assume that $G_0=\mathrm{PSL}_n(q)$ and $H_1\in \mathcal{C}_1\cup\cdots\cup \mathcal{C}_8$. 
Recall Table~\ref{tb:orderofsimplegroups} that 
\[
|\mathrm{Out}(G_0)|=
\begin{cases}
  d\log_p(q)& \text{ if $n=2$},\\
  2d\log_p(q)& \text{ if $n>2$},\\  
\end{cases}
\]
where $d=(n,q-1)$ as Table~\ref{tb:orderofsimplegroups}. 
By Lemma~\ref{lm:q2f},  
\begin{equation}\label{eq:pslout}
\vert \mathrm{Out} (G_0) \vert^2\leq (2d\log_p(q))^2\leq \frac{9}{2}d^{2}q< 5q^3.
\end{equation}

\begin{lemma}[{\cite[Theorem~7 and Proposition 4.7]{AB2015}}]\label{lm:pslc1c8}
Let $G_0=\mathrm{PSL}_n(q)$ and $H_1 \in\mathcal{C}_1 \cup \mathcal{C}_8$. 
Then $H_0$ is large in $G_0$. 
\end{lemma}

\begin{lemma}\label{lm:pslc2}
Let $G_0=\mathrm{PSL}_n(q)$ and $H_1 \in\mathcal{C}_2$.
Then  $H_{1}$ is large in $G_{1}$ if and only if $H_0$ is of type $\mathrm{GL}_{n/t}(q)\wr \mathrm{S}_t$, and  either $t=2$, or $t=3$ and 
\[
q\in \{3, 4, 5, 7, 8, 9, 11, 13, 16, 17, 19, 23, 25, 27, 32, 49, 64, 81, 128\}.
\] 
 In particular,  $H_0$ is large in $G_0$ if and only if $t=2$, or $t=3$ and either $q\in \{5,8,9\}$ and $(n,q-1)=1$, or $q=5$ and $n$ is odd.
\end{lemma}
\begin{proof}
By~\cite[Table 3.5.A]{K-Lie}, $H_1$ is of type $\mathrm{GL}_{m}(q)\wr \mathrm{S}_t$ with $n=mt$ and $t\geq 2$, and $ \mathcal{O}_{1} =\mathrm{Out}(G_0)$ (since $c=1$ and $|\mathcal{O}_{1}|=|\mathrm{Out}(G_0)|/c$).
Moreover, by~\cite[Proposition 2.3.6]{BHRD2013}, we may assume 
\begin{equation}\label{eq:pslc2-qm}
  \text{$q\geq 5$ if $m =1$, and $q \geq 3$ if $m =2$. }
\end{equation} 
By~\cite[Proposition~4.2.9]{K-Lie}, $|H_0|=d^{-1}(q-1)^{t-1}t!|\mathrm{SL}_m(q)^{t}|$.  
Since $ |\mathrm{GL}_m(q)|<(1-q^{-1})q^{m^2}=(q-1)q^{m^2-1}$ (see Lemma~\ref{lm:LUSO}(i)), we have $|\mathrm{SL}_m(q)|<q^{m^2-1}$. Then
 \[|H_0|<d^{-1} q^{m^2t-1}t!.\]
Since $ q^{n^2-2}\leq |G_0|$ (see Lemma~\ref{lm:LUSO1}) and $\vert \mathcal{O}_{1} \vert^2<5d^2q$ (see Eq.~\eqref{eq:pslout}), we conclude from Eq.~\eqref{eq:H1}   that 
\[
q^{m^2t^2-2}<|G_0|\leq |H_0|^3 \vert \mathcal{O}_{1}\vert^2< (d^{-3}q^{3m^2t-3} (t!)^3 ) (5d^{2}q)<5(t!)^3q^{3m^2t-2}  .  
\]
It follows that 
\begin{equation}\label{eq:pslc2qt1}
q^{m^2t(t-3)}< 5(t!)^3.
\end{equation} 
Notice that $t!<2^{t\log_2(\frac{t+1}{2})}$ (see Lemma~\ref{lm:n!}) and  $q^{m^2}\geq 5 =2^{\log_2(5)}$ (see Eq.~\eqref{eq:pslc2-qm}).  
Then we  have
\[
2^{\log_2(5)t(t-3)}\leq q^{m^2t(t-3)}<5(t!)^3<2^{3t\log_2(\frac{t+1}{2})+\log_2(5)},
\]
which implies
\begin{equation}\label{eq:pslc2qt2}
 \log_2(5)t(t-3)<3t\log_2(\frac{t+1}{2})+\log_2(5).
\end{equation} 
A direct verification shows that Eq.~\eqref{eq:pslc2qt2} for integer $t\geq 2$ holds if and only if $t\in \{2,3\}$. 

If $t=2$, then $ H_0  $ is large in $G_0$ by~\cite[Proposition 4.7]{AB2015}.
The case $t=3$ has been dealt with as the example in Subsection~\ref{sec:techniques}.
 \end{proof}  



\begin{lemma}\label{lm:pslc3}
Suppose $G _0=\mathrm{PSL}_n(q)$ and $H_1 \in\mathcal{C}_3$.  
 Then $H_{1}$ is large in $G_{1}$ if and only if  $H_1$ is of type $\mathrm{GL}_{n/r}(q^r)$, where  $r=2$, or $r=3$ with $q \in \{2,3,4,5,7, 8, 9, 11, 16, 27, 32\}$, or $r=5$ and $(G_{1},H_{1})=(\mathrm{PSL}_5(2).2, 31{:}10)$. 
 In particular,  $H_0$ is large in $G_0$ if and only if $r=2$, or $r=3$ and $q\in \{2,3\}$, or $q=5$ and $n$ is odd.
\end{lemma}

\begin{proof}
By~\cite[Table 3.5.A]{K-Lie}, $H_1$ is of type $\mathrm{GL}_m(q^r)$ with $r$ prime and $n=mr$, and $ \mathcal{O}_{1} =\mathrm{Out}(G_0)$. 
Moreover, $|H_0|=d^{-1}(q-1)^{-1}|\mathrm{GL}_m(q^r)|r$ (see~\cite[Proposition 4.3.6]{K-Lie}).
Since $q-1\geq q/2$ and $|\mathrm{GL}_m(q^r)|<q^{m^2r}$, we have $|H_0|\leq  d^{-1}(2rq^{m^2r-1}) $.
Since $ q^{n^2-2}\leq |G|$ and $\vert \mathcal{O}_{1}\vert^2 <9d^2q/2$, we conclude from Eq.~\eqref{eq:H1} that 
\[
q^{m^2r^2-2}<|G_0|\leq  |H_0|^3 \vert\mathcal{O}_{1}\vert^2< (d^{-3} (2r)^3q^{3m^2r-3} ) (9d^2q/2)<36r^3q^{3m^2r-2}  .  
\] 
It follows that
\begin{equation}\label{eq:pslc3qm2s}
r(r-3)<\log_{q^{m^2}}(36r^3).
\end{equation} 
This inequality holds if and only if $r\in \{2,3\}$, or $(q^{m^2},r)=(2, 5) $ (notice that $r$ is a prime).
If $r=2$, then $ H_0  $ is large in $G_0$ by~\cite[Proposition 4.7]{AB2015}. For the case $(q^{m^2},r)=(2, 5) $, we have $(G_{1},H_{1})=(\mathrm{PSL}_5(2).2, 31{:}10)$ (computation shows that $31{:}10$ is indeed large in $\mathrm{PSL}_5(2).2$). 

Assume  $r=3$.
Now Eq.~\eqref{eq:H1} holds if and only if
\[
 \frac{|\mathrm{GL}_m(q^3)|^3 }{|\mathrm{GL}_{3m}(q)|}\geq\frac{(q-1)^2}{3^3\cdot 4\log_p(q)^2}:=h(q)
\] 
By Lemma~\ref{lm:LUSO},
\[
f(q):=\frac{(1-q^{-3}-q^{-6})^3}{(1-q^{-1})(1-q^{-2})}< \frac{|\mathrm{GL}_m(q^3)|^3 }{|\mathrm{GL}_{3m}(q)|} < \frac{(1-q^{-3})^3(1-q^{-6})^3}{ 1-q^{-1}-q^{-2} }:=g(q).
\]
By~\cite[p.202]{AB2015},  $g(q) < 4/3 $ when $q\geq 7$.
Computation shows that  the inequality $h(q)<4/3$ for $q\geq 7$ holds precisely when $q\in \{ 7, 8, 9, 11, 16, 27, 32, 64\}$.
Further computation shows that $h(q)<f(q)$ for each $q\in \{2,3,4,5,7, 8, 9, 11, 16, 27, 32\}$, and $h(64)>g(64)$.
Therefore, in the case $r=3$, the group $H_{1}$ is large in $G_{1}$ if and only if $q\in \{2,3,4,5,7, 8, 9, 11, 16, 27, 32\}$. 
\end{proof}

\begin{lemma}\label{lm:pslc4}
Let $G_0=\mathrm{PSL}_n(q)$ and $H_1\in  \mathcal{C}_4$.
Then $H_{1}$ is not large in $G_{1}$. 
\end{lemma}

\begin{proof}  
By~\cite[Table 3.5.A]{K-Lie}, $H_1$ is of type $\mathrm{GL}_{n_1}(q)\otimes \mathrm{GL}_{n_2}(q) $ with $2\leq  n_1<\sqrt{n} $ and $n=n_1n_2$.
By~\cite[Proposition 4.3.6]{K-Lie},  $|H_0|=d^{-1} c|\mathrm{SL}_{n_1}(q)|\cdot |\mathrm{SL}_{n_2}(q)|$, where $c=(q-1,n_1,n_2)$. 
It follows that
\[
\frac{|H_0|^3|\mathcal{O}_{1}|^2}{|G_0|}=\frac{|\mathrm{SL}_{n_1}(q)|^3\cdot |\mathrm{SL}_{n_2}(q)|^3d^{-3}c^3\cdot (2d\log_p(q)/c)^2 }{d^{-1}|\mathrm{SL}_{n}(q)|}=4\log_p(q)^2c \cdot \frac{|\mathrm{SL}_{n_1}(q)|^3\cdot |\mathrm{SL}_{n_2}(q)|^3}{|\mathrm{SL}_{n}(q)|}.  
\]
By Lemma~\ref{lm:LUSO}, we have 
\[
\frac{1}{2}q^{n^2-1}\leq (1-\frac{q^{-2}}{1-q^{-1}} ) q^{n^2-1}< |\mathrm{SL}_n(q)|<(1-q^{-2})q^{n^2-1}<q^{n^2-1}.
\]
From $|H_0|^3 |\mathcal{O}_1|^2\geq |G_0|$, we have  
\begin{equation}\label{eq:pslc4}
\frac{1}{2} q^{n^2-1}<  |\mathrm{SL}_n(q)|\leq 4\log_p(q)^2c|\mathrm{SL}_{n_1}(q)|^3\cdot |\mathrm{SL}_{n_2}(q)|^3<4\log_p(q)^2c q^{3(n_1^2+n_2^2-2)}).
\end{equation}
Note that $\log_p(q)^2\leq {9}q/{8}$ and $c<q$.
Taking the logarithm base $q$  yields  
\[
n^2-1<\log_{q}(9)+2+3(n_1^2+n_2^2-2).
\]
Recall that $n=n_1n_2$ with $n_2>n_1\geq 2$.
Then 
\[
(n_1^2-3)(n_2^2-3)=n^2-3(n_1^2+n_2^2)+9<6+\log_q(9)<10.
\]
The only possibility is $(n_1,n_2)=(2,3)$.
Now, $c=1$ and $n=6$, and so Eq.~\eqref{eq:pslc4} becomes 
\[
\frac{1}{2}q^{35}<4\log_p(q)^2q^{33} \iff q^{2}<8\log_p(q)^2.
\]
It holds only for $q\in \{2,4,8\}$.
However, a direct computation shows that for $n=6$ and $q\in \{2,4,8\}$, the group $H_1$ is not large in $G_1$.
Hence, the proof is complete. 
\end{proof}
%



\begin{lemma}\label{lm:pslc5}
Let $G_0=\mathrm{PSL}_n(q)$ and $H_1 \in\mathcal{C}_5$.
Then $ H_{1}$ is large in $G_{1}$ if and only if $H_1$ is of type $\mathrm{GL}_n(q^{1/r})$, where $r=2$ or $3$.
Moreover, if $r=2$ then $H_0$ is large in $G_0$; and if $r=3$, then $H_0$ is large in $G_0$ if and only if  $f>1$, where   
\[
f = 
\begin{cases} 
\dfrac{27\,(n, q_0^3 - 1)}{(n, q_0 - 1)^3}, & \text{if } (q_0^2 + q_0 + 1)_3 > 1 \text{ and } (q_0 - 1)_3 \ge n_3 > 1, \\[8pt]
\dfrac{(n, q_0^3 - 1)}{(n, q_0 - 1)^3}, & \text{otherwise.}
\end{cases}
\] 
\end{lemma}

\begin{proof}
According to~\cite[Table 3.5.A]{K-Lie}, $H_1$ is of type $\mathrm{GL}_n(q_0)$, where $q=q_0^r$ for a prime $r$. 
Using~\cite[Proposition 4.5.3]{K-Lie}, $H_0=c(n,q-1)^{-1} |\mathrm{PGL}_n(q_0)|$, where $c= [q_0-1,\frac{q-1}{(q-1,n)}]^{-1}(q-1) $. 
By~\cite[Lemma 1.13.5]{BHRD2013}, $c=(\frac{q-1}{q_0-1},n)$.
Consequently, $H_0\leq |\mathrm{PGL}_n(q_0)|$. 
Since $|\mathrm{PGL}_n(q_0)|=(1-q_0)^{-1} |\mathrm{GL}_n(q_0)|$ and $|\mathrm{GL}_n(q_0)|\leq (1-q_0^{-1})(1-q_0^{-2})q_0^{n^2 }$ (see Lemma~\ref{lm:LUSO}), we have  $|H_0|\leq q_0^{n^2-1}$. 
Notice that  $|G_0|\geq q^{n^2-2}$ and  $|\mathcal{O}_{1}|^2 < 5q^3$.  
From Eq.~\eqref{eq:H1} we have 
\[
q_0^{r(n^2-2)}= q^{n^2-2}<|G_0|\leq  |H_0|^3 |\mathcal{O}_{1}|^2< q_0^{3(n^2-1)}(4d^2\log_p(q)^2)<5 q_0^{3n^2 }  .  
\]  
It follows that
\begin{equation}\label{eq:pslc5}
n^2(r-3)-2r <\log_{q_0}(5).
\end{equation} 
Notice that  $r$ is a prime and $r\mid n$.
It is clear that $n^2(r-3)-2r $ is increasing on $r$.
If $r\geq 5$, then $n^2(r-3)-2r-3>5^2\cdot 2-10-3=37$, while $\log_{q_0}(5)\leq \log_2(5)<3$, a contradiction.
Thus, Eq.~\eqref{eq:pslc5} holds if and only if $r\in \{2,3\}$.

We have proved that if $H_1$ is large in $G_1$ then  $r\in \{2,3\}$.
According to~\cite[Proposition 4.7]{AB2015}, if $r=2$, then $H_0$ is large in $G_0$, and if $r=3$, then $H_0$ is large in $G_0$ if and only if  $f>1$. 
Thus, it suffices to show that $H_1$ is large in $G_1$ when $r=3$. 

Assume $r=3$. Then 
\begin{align*}
  \frac{|H_0|^3|\mathcal{O}_{1}|^2}{|G_0|}
  &=\frac{\frac{c^3|\mathrm{PGL}_n(q_0)|^3}{(n,q -1)^3}\cdot  \frac{4\log_p(q)^2(n,q-1)^2}{c^2}}{\frac{1}{(n,q-1)}|\mathrm{PGL}_n(q_0^3)|}
  = 4\log_p(q)^2c    \cdot \frac{|\mathrm{PGL}_n(q_0)|^3}{ |\mathrm{PGL}_n(q)|}\\
  & = 4\log_p(q)^2 \left(\frac{q-1}{q_0-1},n\right) \cdot \frac{q-1}{(q_0-1)^3}\cdot \frac{|\mathrm{GL}_n(q_0)|^3}{ |\mathrm{GL}_n(q)|}\\
  &= 36\log_p(q_0)^2 \left(\frac{q_0^3-1}{q_0-1},n\right) \cdot \frac{q_0^3-1}{(q_0-1)^3}\cdot \frac{|\mathrm{GL}_n(q_0)|^3}{ |\mathrm{GL}_n(q_0^3)|}\\
  &>36\cdot\frac{q_0^3-1}{(q_0-1)^3}\cdot \frac{|\mathrm{GL}_n(q_0)|^3}{ |\mathrm{GL}_n(q_0^3)|}\\
  &>36\cdot \frac{q_0^3-1}{(q_0-1)^3}\cdot \frac{(1-q_0^{-1}-q_0^{-2})^3}{(1-q_0^{-3})(1-q_0^{-6})}:=h(q_0).
\end{align*} 
Note that the last inequality follows from the upper bound for $|\mathrm{GL}_n(q_0)|$ and the lower bound for $|\mathrm{GL}_n(q_0^3)|$ given in Lemma~\ref{lm:LUSO}. 
One can easily verify that $h(q_0)$ is increasing in $q_0$, and $h(2)=32/7>1$.
Therefore, $|H_0|^3|\mathcal{O}_{1}|^2>|G_0|$ always holds; that is, $H_1$ is always large in $G_0$ when $r=3$, as expected.
\end{proof}

\begin{lemma}\label{lm:pslc6}
Let $G_0=\mathrm{PSL}_n(q)$ and $H_1 \in\mathcal{C}_6$.
Then $H_{1}$ is large in $G_{1}$ if and only if  $(G_0,H_0)$ is one of the following
\[
\begin{split}
& (\mathrm{PSL}_{4}(5),2^{4}.\mathrm{A}_{6} ), (\mathrm{PSL}_{3}(7),3^{2}.\mathrm{Q}_{8} ),\\ 
& (\mathrm{PSL}_{2}(q),\mathrm{S}_{4}) \text{ with $q\in \{7,17,23,31,41,47\}$},\\
 &(\mathrm{PSL}_{2}(q),\mathrm{A}_{4}) \text{ with $q\in \{5,11,13,19,37\}$}.
 \end{split}
\]   
In particular, $H_0$ is large in $G_0$ if and only of $(G_0,H_0)=(\mathrm{PSL}_{4}(5),2^{4}.\mathrm{A}_{6} )$, or $  (\mathrm{PSL}_{2}(q),\mathrm{S}_{4})$ with $q\in \{7,17,23 \}$, or $  (\mathrm{PSL}_{2}(q),\mathrm{A}_{4})$ with $q\in \{5,11,13 \}$.
\end{lemma} 

\begin{proof}   
By~\cite[Table 4.6.B]{K-Lie},  $H_1$ is one of the following types:
\begin{enumerate}[\rm (i)]
  \item $r^{1+2m}.\mathrm{Sp}_{2m}(r)$, where $n=r^{m}$, $m\geq 1$, and $r$ is an odd prime divisor of $q-1$, and $er$ is odd ($e=\log_p(q)$ is the smallest integer such that $r\mid (q^e-1)$);
  \item $(4\circ 2^{1+2m}).\mathrm{Sp}_{2m}(2)$, where $n=2^{m}\geq 4$, $q=p$ is odd;
  \item $2_{-}^{1+2}.\mathrm{O}^{-}_{2}(2)$, where $ n=2$ and $q=p$ is odd.
\end{enumerate}

Assume that $H_1$ is of type  $r^{1+2m}.\mathrm{Sp}_{2m}(r)$.
By~\cite[Proposition 4.6.5]{K-Lie}, $|H_0|\leq r^{2m}|\mathrm{Sp}_{2m}(r)|$.
Since $|\mathrm{Sp}_{2m}(r)|<r^{m(2m-1)}$ (see Lemma~\ref{lm:LUSO}), we derive that $|H_0|\leq  r^{2m+m(2m+1)}=r^{m(2m+3)} $. 
By Eq.~\eqref{eq:H1} we have 
\[
  q^{n^2-2}<|G_0|\leq  |H_0|^3 \vert \mathcal{O}_{1}\vert^2< r^{3m(2m+3)}(5q^3) .  
\] 
It follows that
\begin{equation}\label{eq:pslc61}
r^{2m}-5=n^2-5 <3m(2m+3)\log_q(r)+\log_{q}(5).
\end{equation} 
Since $r\geq 3$ and $r \mid q-1$, we have $q\geq 4$ and $\log_q(r)<1$.
Hence, 
\begin{equation}\label{eq:pslc62}
 r^{2m}-3m(2m+3)-7<0.
\end{equation} 
If $m\geq 2$, then by the binomial expansion, we have
\begin{align*}
  r^{2m}=(r-2+2)^{2m}&>2^{2m}+2m\cdot 2^{2m-1}(r-1)+ m(2m-1)\cdot 2^{2m-2}(r-2)^2\\
  &> 16+16m+4m(2m-1)> 7+3m(2m+3),  
\end{align*}
which is a contradiction.
Therefore, $m=1$.
From Eq.~\eqref{eq:pslc62} we have $r^2-22<0$, forcing $r=3$. 
Now, $G_0=\mathrm{PSL}_3(q)$.
By~\cite[Table~8.3]{BHRD2013}, $H_0=3^{2}{:}\mathrm{Q}_8$, $|\mathcal{O}_{1}|=2d$ and $q=p \equiv 1\pmod{3}$. 
One can verify that only $q=7$ satisfying Eq.~\eqref{eq:H1}. 
Hence, we obtain $(G_0,H_0)=(\mathrm{PSL}_3(7),3^{2}{:}\mathrm{Q}_8)$ in the case where $H_1$ is of type  $r^{1+2m}.\mathrm{Sp}_{2m}(r)$. Note that $|\mathrm{PSL}_3(7)|/|3^{2}{:}\mathrm{Q}_8|^3\approx 5$ and so $3^{2}{:}\mathrm{Q}_8$ is not large in $\mathrm{PSL}_3(7)$.
 
Assume that $H_1$ is of type   $ (4\circ 2^{1+2m}).\mathrm{Sp}_{2m}(2) $.
By~\cite[Proposition 4.6.6]{K-Lie}, $|H_0|\leq 2^{2m}|\mathrm{Sp}_{2m}(2)|$.
Then Eq~\eqref{eq:pslc61} also holds here.
Now $q\geq 3$ and $m\geq 2$.
Thus,
\[
2^{2m}-5<m(6m+9)+\log_3(5).
\]
It holds if and only if $m\in \{2,3\}$.

Let $m=2$. Then $n=4$. By~\cite[Table~8.8]{BHRD2013}, either  
\begin{itemize}
  \item $|H_0|=2^4\cdot |\mathrm{S}_6|$, $q=p=1\pmod{8}$ and $c=4$; or
  \item $|H_0|=2^4\cdot |\mathrm{A}_6|$, $q=p=5\pmod{8}$ and $c=2$.
\end{itemize}
Direct verification shows that  Eq.~\eqref{eq:H1} holds if and only if $q=5$, giving $(G_0,H_0)=(\mathrm{PSL}_4(5),2^4.\mathrm{A}_6)$.

Let $m=3$. Then $n=8$, $H_0=2^6.\mathrm{Sp}_6(2)$ and $q=p=1\pmod{4}$, and $c=d=(8,q-1)$. Direct verification shows that  Eq.~\eqref{eq:H1} does not hold.

Assume that $H_1$ is of type  $2_{-}^{1+2}.\mathrm{O}^{-}_{2}(2)$. Now $n=2$.
By~\cite[Table~8.1]{BHRD2013}, either
\begin{itemize}
  \item $|H_0|= |\mathrm{S}_4|$, $q=p=\pm 1\pmod{8}$ and $c=2$; or
  \item $|H_0|= |\mathrm{A}_4|$, $q=p=\pm 3,5,\pm 11,\pm 13,\pm 19\pmod{40}$ and $c=1$.
\end{itemize} 
Direct verification gives the remaining pairs $(G_0,H_0)$ in the lemma. 
\end{proof}

\begin{lemma}
Let $G_0=\mathrm{PSL}_n(q)$ and $H_1 \in \mathcal{C}_7 $.
Then $H_{1}$ is not large in $G_{1}$. 
\end{lemma}

\begin{proof}
Suppose for a contradiction that $H_{1}$ is large in $G_{1}$. 
By~\cite[Table 3.5.A and Proposition 4.7.3]{K-Lie}, $H_1$ is of type $\mathrm{GL}_{m}(q) \wr \mathrm{S}_t$, where $n=m^t$ with $m\geq 3$ and $t \geq 2$, and $|H_0|<|\mathrm{SL}_{m}(q)|^{t}t!$.
Since $|\mathrm{SL}_{m}(q)|<q^{m^2-1}$, $|G_0|\geq q^{n^2-2}$, $\vert \mathcal{O}_{1} \vert^2<5q^3$ and $t!<2^{t\log_2(\frac{t+1}{2})}$, we conclude from Eq.~\eqref{eq:H1} that 
\[
q^{m^{2t}-2}<|G_0|\leq  |H_0|^3 \vert \mathcal{O}_{1} \vert^2< (q^{t(m^2-1)}t!)^3\cdot (5q^3)<5q^{3t(m^2-1)+3}\cdot 2^{3t\log_2(\frac{t+1}{2})} .  
\] 
It follows that 
\[
 m^{2t}-3t(m^2-1)-5<\log_2(5) +3t\log_2\left(\frac{t+1}{2}\right).
\]
It is clear that $ m^{2t}-3t(m^2-1) $ is increasing in $m$. Since $m\geq 3$, we have
\[
 3^{2t}-24t-5 <\log_2(5) +3t\log_2\left(\frac{t+1}{2}\right). 
\]
However, there is no $t\geq 2$ satisfying the above inequality. 
This contradiction proves the lemma.
\end{proof}

\subsection{Geometric  subgroups of unitary groups} \label{sec:proof-psu}
In this subsection, we assume that $G_0=\mathrm{PSU}_n(q)$ with $n\geq 3$.  
Now $|\mathrm{Out}(G_0)|= 2d\log_p(q)$, 
where $d=(n,q+1)$. By~Lemma~\ref{lm:q2f}, we conclude
\begin{equation}\label{eq:psuout-1}
  \vert \mathcal{O}_{1}\vert^2\leq \vert \mathrm{Out} (G_0)\vert^2 \leq  (2d\log_p(q))^2 <  \frac{9}{2}d^2q  .
  \end{equation}  
Notice that $2d=2(n,q+1)\leq 2(q+1)\leq 3q$.  
If $q\neq 8$, then $(\log_p(q))^2\leq q$ by Lemma~\ref{lm:q2f}, and hence $\vert \mathrm{Out} (G_0) \vert^2\leq 9q^3$. 
If $q=8$, then $(\log_p(q))^2=9$  and so $\vert \mathrm{Out} (G_0) \vert^2 \leq 2^2\cdot 9^2\cdot 9=2916<4608=9q^3$.  
Therefore, $\vert \mathrm{Out} (G_0)\vert^2<9q^3$. Hence, 
\begin{equation}\label{eq:psuout}
\vert \mathcal{O}_{1}\vert^2\leq \vert \mathrm{Out} (G_0)\vert^2 < \text{min}\{\frac{9}{2}d^2q, 9q^3\}.
\end{equation}

\begin{lemma}[{\cite[Proposition 4.17]{AB2015}}]\label{lm:psuc1}
Let $G_0=\mathrm{PSU}_n(q)$ and $H_1 \in\mathcal{C}_1 $. 
Then $H_0$ is large in $G_0$. 
\end{lemma}

\begin{lemma}\label{lm:psuc2}
Let $G_0=\mathrm{PSU}_n(q)$ and $H_1 \in\mathcal{C}_2$. 
Then $H_{1}$ is large in $G_{1}$ if and only if one of the following holds:
\begin{enumerate}[\rm (a)]
  \item $H_1$ is of type $\mathrm{GL}_{n/2}(q^2).2$;
   \item  $H_1$ is of type $\mathrm{GU}_{n/t}(q)\wr \mathrm{S}_t$, and  one of the following holds:
  \begin{enumerate}[\rm (b.1)]  
  \item $t=2$;
  \item $t=3$ and $q \in \{2, 3, 4, 5, 7, 8, 9, 11, 13, 16, 17, 19, 23, 25, 27, 29, 32, 49, 64, 81, 128\}$;
  \item $4\leq t\leq 11$, and $(q,m,t)$ is as follows:
  \[
   \begin{split}
  &   ( 2, 1, 4 ),
      ( 2, 1, 5 ),
      ( 2, 1, 6 ),
      ( 2, 1, 7 ),
      ( 2, 1, 8 ),
      ( 2, 1, 9 ),
      ( 2, 1, 10 ),
      ( 2, 1, 11 ), \\&
      ( 3, 1, 4 ), 
      ( 3, 1, 5 ),
      ( 3, 1, 6 ),
      ( 4, 1, 4 ),
      ( 4, 1, 5 ),
      ( 5, 1, 4 ),
      ( 7, 1, 4 ),
      ( 8, 1, 4 ),
      ( 9, 1, 4 ).
  \end{split}
  \]
\end{enumerate}  
\end{enumerate}  
In particular, $H_0$ is large in $G_0$ if and only if 
 $H_1$ is of type $\mathrm{GL}_{n/2}(q^2).2$, or of type $\mathrm{GU}_{n/t}(q)\wr \mathrm{S}_t$,  where one of the following holds:
\begin{enumerate}[\rm (i)]
  \item $t=2$;
  \item  $t=3$ and either $q\in \{2,3,4\}$, or $(q,(n,q+1))=(5,3),(7,1),(7,2)$, $(9,1)$, $(9,2)$, $(13,1)$ or $(16,1)$;
  \item $4\leq t=n\leq 11$, and either $q=2$, or $(n,q)=(6,3)$, $(5,3)$, $(4,3)$, $(4,4)$ or $(4,5)$. 
\end{enumerate} 
\end{lemma}
\begin{proof}
According to~\cite[Table 3.5.B]{K-Lie}, $H_1$ is of  type $\mathrm{GL}_{n/2}(q^2).2$ or $\mathrm{GU}_m(q)\wr \mathrm{S}_t$.  
If $H_1$ is of type $\mathrm{GL}_{n/2}(q^2).2$, then $H_0$ is large in $G_0$, see~\cite[Proposition 4.17]{AB2015}. 

Assume that  $H_1$ is of type $\mathrm{GU}_m(q)\wr \mathrm{S}_t$, where $n=mt$ and $t\geq 2$. 
By~\cite[Proposition 4.3.6]{K-Lie},  $ \mathcal{O}_{1} =\mathrm{Out}(G_0)$, $|H_0|=d^{-1}(q+1)^{-1}|\mathrm{GU}_m(q)|^t t!$.   
Notice that $  |G_0|>  (1-q^{-1})q^{n^2-2}>q^{n^2-2}/2$ (see Lemma~\ref{lm:LUSO1}) and $\vert \mathcal{O}_{1} \vert^2\leq 9d^2q/2 $ (see Eq.~\eqref{eq:psuout}). 
By Lemma~\ref{lm:LUSO}(i),
\[
|\mathrm{GU}_m(q)|\leq (1+q^{-1} ) (1-q^{-2} ) (1+q^{-3} )q^{m^{2}}\leq (1+q^{-1} ) q^{m^2}.
\]
 Since $(1+q^{-1})=q^{-1}(q+q^{-1})$, we have 
 \[|H_0|<d^{-1} (1+q^{-1})^{t-1}q^{m^2t-1}t!.\] 
Then from Eq.~\eqref{eq:H1} we have  that 
\[
\begin{split}
 \frac{1}{2}q^{m^2t^2-2}  <|G_0|\leq |H_0|^3 \vert \mathcal{O}_{1}\vert^2 &< (d^{-3}(1+q^{-1})^{3(t-1)} q^{3m^2t-3} (t!)^3 ) (\frac{9}{2}d^{2}q)\\
&<9(t!)^3(1+q^{-1})^{3(t-1)}q^{3m^2t-2} .
\end{split}  
\]
It follows that 
\begin{equation}\label{eq:psuc2-1}
  q^{m^2t(t-3)}< 9(1+q^{-1})^{3(t-1)}(t!)^3.
\end{equation} 
Clearly, the above inequality holds when $t\in \{2,3\}$.   
If $t=2$, then $ H_0  $ is large in $G_0$ by~\cite[Proposition 4.17]{AB2015}.

Assume $t=3$. Notice that  $|G_0|=d^{-1}(q+1)^{-1}|\mathrm{GU}_{3m}(q)|$. 
Hence,
\[
\frac{|H_0|^3|\mathcal{O}_{1}|^2}{|G_0|}=\frac{(3!)^3\cdot (4d^2\log_p(q)^2)}{d^2(q+1)^2} \cdot \frac{|\mathrm{GU}_{ m}(q)|^9}{|\mathrm{GU}_{3m}(q)|}=\frac{864\log_p(q)^2}{(q+1)^2} \cdot \frac{|\mathrm{GU}_{ m}(q)|^9}{|\mathrm{GU}_{3m}(q)|}.
\]
Then  Eq.~\eqref{eq:H1} is equivalent to 
\[ \frac{|\mathrm{GU}_{ m}(q)|^9}{|\mathrm{GU}_{3m}(q)|}\geq \frac{(q+1)^2}{864\log_p(q)^2}:=h(q).
\] 
By Lemma~\ref{lm:LUSO}, $(1+q^{-1} ) (1-q^{-2} )q^{a^{2}} <|\mathrm{GU}_{a}(q) |<(1+q^{-1} )  q^{a^{2}}$, hence 
\[
f(q):=\frac{(1+q^{-1})^9(1-q^{-2})^9}{ 1+q^{-1}  }<\frac{|\mathrm{GU}_{ m}(q)|^9}{|\mathrm{GU}_{3m}(q)|}<\frac{( 1+q^{-1} )^9 }{(1+q^{-1}) (1-q^{-2})}=:g(q).
\] 
If $q\geq 200$, then $h(q)>(q+1)^2/864>46$, while $g(q)<4 (1+200^{-1})^8<5$, implying that Eq.~\eqref{eq:H1} cannot happen. 
Therefore, $q<200$. Computation shows that  $h(q)<g(q)$ holds if and only if 
\[
q \in \{2, 3, 4, 5, 7, 8, 9, 11, 13, 16, 17, 19, 23, 25, 27, 29, 32, 49, 64, 81, 128\}.
\]
Moreover, for those $q$ above, $h(q)<f(q)$ also holds,  implying that $H_1$ is indeed large in $G_1$.  This proves part (b.2).

Assume $t\geq 4$. 
Set 
\[
f(q,m,t):=t(t-3)-\log_{q^{m^2}}(9)-3(t-1)\log_{q^{m^2}}(1+q^{-1})-3\log_{q^{m^2}}(t!).
\] 
Then Eq.~\eqref{eq:psuc2-1} is equivalent to $f(q,m,t)<0$.
It is clear that $f(q,m,t)$ is increasing in both $q$ and $m$.
Therefore,  
\[
  f(2,1,t)\leq f(2,m,t)\leq f(q,m,t)<0.
\]
Note that \( f(2,1,t) \) is a univariate function:  
\[
f(2,1,t) = t(t-3) - \log_{2}(9) - 3(t-1)\log_{2}\left(1+2^{-1}\right) - 3\log_{2}(t!).
\]  
Since $2\log_{2}\left(1+2^{-1}\right) < 1.8$ and  $t! < 2^{t\log_2\left(\frac{t+1}{2}\right)}$,
we observe that if \( f(2,1,t) < 0 \), then \( g(t) < 0 \) must necessarily hold, where  
\[
g(t) := t(t-3) - 1.8(t-1) - 3t\log_2\left(\frac{t+1}{2}\right) - \log_{2}(9).
\]  
It is straightforward to verify that, for sufficiently large \( t \) (e.g., \( t > 100 \)), we have $\frac{t+1 }{2} \geq 3 \log_2\left(\frac{t+1}{2}\right)$, and  further $\frac{t(t-7)}{2} > 1.8(t-1) + \log_{2}(9)$.  
Consequently, the solution set for \( g(t) < 0 \) is confined to \( 4 \leq t \leq 100 \).
Direct verification shows that  \( f(2,1,t) < 0 \) with $4 \leq t \leq 100$ holds  if and only if  $4\leq t\leq 12$. Then for each $4\leq t\leq 12$, we compute the values for  $m$ such that $f(2,m,t)<0$.
It shows that  either $(m,t)=(2,4)$, or $m=1$ and $4\leq t\leq 12$. 
Then for each pair $(m,t)$, we compute those $q$ such that $f(q,m,t)<0$. 
Computation shows that $(q,m,t)$ is one of the following:  $( 2 ,  2 ,  4 )$, $(2,1,i)$ with $4\leq i\leq 12$, $(j,1,4)$ with $3\leq j \leq 19$, $(k,1,5)$ with $3\leq k \leq 5$, $( 3 ,  1 ,  6)$, $( 3 ,  1 ,  7)$, $( 4 ,  1 ,  6)$.
For those triples $(q,m,t)$, direct computation shows that Eq.~\eqref{eq:H1} holds if and only if $(q,m,t)$ is as follows:
\[
 \begin{split}
&   ( 2, 1, 4 ),
    ( 2, 1, 5 ),
    ( 2, 1, 6 ),
    ( 2, 1, 7 ),
    ( 2, 1, 8 ),
    ( 2, 1, 9 ),
    ( 2, 1, 10 ),
    ( 2, 1, 11 ), \\&
    ( 3, 1, 4 ),
    ( 3, 1, 5 ),
    ( 3, 1, 6 ),
    ( 4, 1, 4 ),
    ( 4, 1, 5 ),
    ( 5, 1, 4 ),
    ( 7, 1, 4 ),
    ( 8, 1, 4 ),
    ( 9, 1, 4 ).
\end{split}
\]
This proves part (b.3).  
\end{proof}

\begin{lemma}\label{lm:psuc3}
Suppose $G _0=\mathrm{PSU}_n(q)$ and $H_1 \in\mathcal{C}_3$.
Then $H_{1}$ is large in $G_{1}$ if and only if $H_0$ is of type $\mathrm{GU}_{n/r}(q^r)$ with $r=3$ and $q \in \{2, 3, 4, 5, 7, 8, 9, 16, 27, 32\}$. 
In particular, $H_0$ is large in $G_0$ if and only if $q=3$ and $n$ is odd.
\end{lemma}

\begin{proof}
By~\cite[Table 3.5.B]{K-Lie}, $H_1$ is of type $\mathrm{GU}_{m}(q^r)$, where  $r$ is an odd prime and $n=mr$, and $ \mathcal{O}_{1} =\mathrm{Out}(G_0)$. By~\cite[Proposition 4.3.6]{K-Lie},  $|H_0|=d^{-1}(q+1)^{-1}|\mathrm{GU}_m(q^r)|r $.
Since $|\mathrm{GU}_m(q^r)|<(1+q^{-1})q^{m^2r}$, we have $|H_0|\leq  d^{-1}r q^{m^2r-1}  $. 
Since $  |G_0| > (1-q^{-1})q^{n^2-2}\geq q^{n^2-2}/2$ and $\vert \mathcal{O}_{1}\vert^2 <9d^2q/2$, we conclude from Eq.~\eqref{eq:H1} that 
\[
\frac{1}{2}q^{m^2r^2-2} <|G_0|\leq |H_0|^3 \vert\mathcal{O}_{1}\vert^2< (d^{-3} r^3q^{3m^2r-3} ) (\frac{9}{2}d^2q )<\frac{9}{2}t^3q^{3m^2r-2} .  
\] 
It follows that
\begin{equation} \label{eq:psuc3qm2s}
r(r-3)<\log_{q^{m^2}}(9r^3).
\end{equation} 
This inequality holds if and only if $r=3$, or $(q^{m^2},r)=(2, 5) $ (notice that $r$ is an odd prime).
 For the case $(q^{m^2},r)=(2, 5) $, we have $|G_0|=|\mathrm{PSU}_5(2)|=13685760$, $|H_0|^3=(|11{:}5|)^3 =166375$ and $|\mathcal{O}_{1}|=2$, for which Eq.~\eqref{eq:H1} does not hold. 
  
Assume  $r=3$.
Now Eq.~\eqref{eq:H1} holds if and only if
\[
 \frac{|\mathrm{GU}_m(q^3)|^3 }{|\mathrm{GU}_{3m}(q)|}>\frac{(q+1)^2}{3^3\cdot 4\log_p(q)^2}:=h(q)
\] 
By Lemma~\ref{lm:LUSO}, $(1+q^{-1} ) (1-q^{-2} )q^{a^{2}} <|\mathrm{GU}_{a}(q) |<(1+q^{-1} )  q^{a^{2}}$, hence  
\[
f(q):=\frac{(1+q^{-3})^3(1-q^{-6})^3}{ 1+q^{-1}  }<\frac{|\mathrm{GU}_{ m}(q^3)|^3}{|\mathrm{GU}_{3m}(q)|}<\frac{( 1+q^{-3} )^3 }{(1+q^{-1}) (1-q^{-2})}=:g(q).
\] 
If $q\geq 1080$, then $h(q)>(q+1)^2/108q>10$, while $g(q)<4\cdot (1+200^{-3})^3<5$, implying that Eq.~\eqref{eq:H1} cannot happen. 
Therefore, $q<200$. Computation shows that  $h(q)<g(q)$ holds if and only if 
\[
q \in \{2, 3, 4, 5, 7, 8, 9, 16, 27, 32\}.
\]
Moreover, for those $q$ above, $h(q)<f(q)$ also holds, implying that $H_1$ is indeed large in $G_1$.   
\end{proof}

\begin{lemma}\label{lm:psuc4}
Let $G_0=\mathrm{PSU}_n(q)$ and $H_1\in  \mathcal{C}_4$.
Then $H_{1}$ is not large in $G_{1}$. 
\end{lemma}

\begin{proof}
Suppose for a contradiction that $H_{1}$ is large in $G_{1}$.
By~\cite[Table 3.5.B and Proposition 4.4.10]{K-Lie}, $H_1$ is of type $\mathrm{GU}_{n_1}(q)\otimes \mathrm{GU}_{n_2}(q) $ with $1< {n_1}<\sqrt{n} $ and $n={n_1}{n_2}$, and  furthermore, $|H_0|=d^{-1}c|\mathrm{SU}_{n_1}(q)|\, | \mathrm{SU}_{n_2}(q)|$, where $c=(q+1,{n_1},{n_2})$. 
By Lemma~\ref{lm:LUSO}, we conclude  $|\mathrm{SU}_a(q)|\leq (1-q^{-2})(1+q^{-3})q^{a^2-1}< q^{a^2-1}$.
Notice that $|G_0|> q^{n^2-2}/2$ and $\vert \mathcal{O}_{1}\vert^2\leq   4d^2\log_p(q)^2$.
Since $1<n_1<\sqrt{n}$, we have $n\geq 6$ and $n_2\geq 3$.

Assume $(n_1,n_2)=(2,3)$. Then $|H_0|=d^{-1}|\mathrm{SU}_{n_1}(q)|\cdot |\mathrm{SU}_{n_2}(q)|\leq d^{-1}q^{11}$.
From Eq.~\eqref{eq:H1} we have
\[
\frac{1}{2}q^{34} <|G_0|\leq  |H_0|^3 \vert \mathcal{O}_{1}  \vert^2< (d^{-3}  q^{33} ) (4d^2\log_p(q)^2) < 4 \log_p(q)^2 q^{33 }  .
\]  
It follows $q<8\log_p(q)^2$, which holds only for 
\[
q\in \{ 2, 3, 4, 5, 7, 8, 9, 16, 25, 27, 32, 64, 81, 128, 256, 512\}.
\]
Computation shows that Eq.~\eqref{eq:H1} does not hold for any $q$ above.
Therefore, the case  $(n_1,n_2)=(2,3)$ cannot happen.

Assume that $(n_1,n_2)\neq (2,3)$. 
Now
\[
| H_0| \leq \frac{|\mathrm{SU}_{n_1}(q)||\mathrm{SU}_{n_2}(q)| (q+1) }{d }\leq \frac{q^{{n_1}^2+{n_2}^2-2}}{ d} \cdot \frac{3q}{2}=\frac{3q^{{n_1}^2+{n_2}^2-1}}{2d}.
\] 
Since  $\vert \mathcal{O}_{1}\vert^2\leq \vert \mathrm{Out} (G_0) \vert^2 \leq 9d^2q/2$, we conclude from Eq.~\eqref{eq:H1} that 
\[
  \frac{1}{2} q^{{n_1}^2{n_2}^2-2} <|G_0|\leq  |H_0|^3 \vert \mathcal{O}_{1}  \vert^2< ( (\frac{3}{2d})^3 q^{3({n_1}^2+{n_2}^2-1)} ) (\frac{9}{2}d^2q ) < \frac{243 }{16}q^{3({n_1}m^2+{n_2}^2-2) }  . 
\] 
It follows that 
\begin{equation}\label{eq:psuc4}
  {n_1}^2{n_2}^2-3({n_1}^2+{n_2}^2) < \log_q( \frac{243}{8})<\log_q( 31).
\end{equation} 
Since $({n_1},{n_2})\neq (2,3)$ and $1<{n_1}<\sqrt{{n_1}{n_2}}$, we have ${n_2}\geq 4$.
If ${n_2}\geq 5$, then 
\[ {n_1}^2{n_2}^2-3({n_1}^2+{n_2}^2)=({n_1}^2-3){n_2}^2-3{n_1}^2\geq 22{n_1}^2-75\geq 22\cdot 4-75=13>\log_2(31)\geq \log_q(31),\] a contradiction.
Therefore, ${n_2}=4$. Now ${n_1}^2{n_2}^2-3({n_1}^2+{n_2}^2)  =13{n_1}^2-48$.
It implies that Eq.~\eqref{eq:psuc4} holds only for $({n_1},{n_2},q)=(2,4,2)$. 
However, by~\cite[Tables~8.46]{BHRD2013}, when $q=2$, such maximal subgroup $H_{1}$ of type  $\mathrm{GU}_2(q)\otimes \mathrm{GU}_4(q) $  do not exist.   
\end{proof}

\begin{lemma}\label{lm:psuc5}
Let $G_0=\mathrm{PSU}_n(q)$ and $H_1 \in\mathcal{C}_5$. 
If $H_1$ is large in $G_1$, then $H_0$ is of type $\mathrm{Sp}_n(q)$, or $\mathrm{O}_{n}^{\epsilon}(q)$  with $\epsilon\in \{+,-,\circ\}$ and $q$ odd, or $\mathrm{GU}_n(q^{1/r})$, where $r=2$ or $3$.  
In particular, if $G_0$ of type $\mathrm{Sp}_n(q)$, or $\mathrm{O}_{n}^{\epsilon}(q)$, or $\mathrm{GU}_n(q^{1/r})$ with $r=2$, then $H_0$ is large in $G_0$; if $G_0$ of type  $\mathrm{GU}_n(q^{1/r})$ with $r=3$, then $H_0$ is large in $G_0$ if and only if $f>1$, where   
\[
f = 
\begin{cases} 
\dfrac{27\,(n, q_0^3 + 1)}{(n, q_0 + 1)^3}, & \text{if } (q_0^2 - q_0 + 1)_3 > 1 \text{ and } (q_0 + 1)_3 \ge n_3 > 1, \\[8pt]
\dfrac{(n, q_0^3 + 1)}{(n, q_0 + 1)^3}, & \text{otherwise.}
\end{cases}
\]
\end{lemma}

\begin{proof}
According to \cite[Table 3.5.B]{K-Lie}, a \(\mathcal{C}_5\)-subgroup of \(G_0\) can only be of one of the following types: \(\mathrm{Sp}_n(q)\), \(\mathrm{O}_{n}^{\epsilon}(q)\) with $\epsilon\in \{+,-,\circ\}$ and $q$ odd, or  \(\mathrm{GU}_n(q_0)\) with \(q = q_0^r\) for a prime \(r\).
The ``in particular" part of the statement is taken from \cite[Proposition~4.17]{AB2015}.
 
Assume that  $H_1$ is of type $\mathrm{GU}_n(q_0)$ with $q=q_0^r$ and $r \geq 3$. 
According to~\cite[Proposition 4.5.3]{K-Lie}, $|H_0|=c(n,q+1)^{-1} |\mathrm{PGU}_n(q_0)|$, where $c= [q_0+1,\frac{q+1}{(q+1,n)}]^{-1}(q+1) $. 
Again by~\cite[Lemma 1.13.5]{BHRD2013}, $c=(\frac{q+1}{q_0+1},n)$.
Hence, $|H_0|\leq |\mathrm{PGU}_n(q_0)|$.
Since $|\mathrm{PGU}_n(q_0)|=|\mathrm{SU}_n(q_0)|=(1+q_0)^{-1} |\mathrm{GU}_n(q_0)|$ and $|\mathrm{GU}_n(q_0)|\leq (1+q_0^{-1})(1-q_0^{-2})q_0^{n^2}$ (see Lemma~\ref{lm:LUSO}), we have  $|H_0|\leq q_0^{n^2-1}$.  
Notice that $|G_0|>(1-q^{-1})q^{n^2-1}\geq \frac{7}{8}q^{n^2-1}$ (as $q=q_0^r\geq 8$).
From Eq.~\eqref{eq:H1} we have  
\[
\frac{7}{8} q_0^{r(n^2-2)}= \frac{7}{8}q^{n^2-2}<|G_0|\leq |H_0|^3 \vert \mathcal{O}_{1} \vert^2< q_0^{3n^2-3}(5q^3)<5 q_0^{3n^2 }  .  
\]  
Taking logarithms base \(q_0\) yields
\begin{equation}\label{eq:psuc5}
n^2(r-3)-2r <\log_{q_0}(\frac{40}{7})<\log_{q_0}(6).
\end{equation} 
Notice that  $r$ is a prime and $r\mid n$.
If $r\geq 5$, then $n^2(r-3)-2r\geq r^2(r-3)-2r \geq 5^2\cdot 2-10 >\log_{2}(6)>\log_{q_0}(6)$, a contradiction. 
Thus, if \(H_1\) is large in \(G_1\), the  parameter \(r\) can only be \(3\). 

Now, let $r=3$. Then 
\begin{align*}
  \frac{|H_0|^3|\mathcal{O}_{1}|^2}{|G_0|}
  &=\frac{\frac{c^3|\mathrm{PGU}_n(q_0)|^3}{(n,q+1)^3}\cdot  \frac{4\log_p(q)^2(n,q+1)^2}{c^2}}{\frac{1}{(n,q+1)}|\mathrm{PGU}_n(q_0^3)|}
  = 4\log_p(q)^2c    \cdot \frac{|\mathrm{PGU}_n(q_0)|^3}{ |\mathrm{PGU}_n(q)|}\\
  & = 4\log_p(q)^2 \left(\frac{q+1}{q_0+1},n\right) \cdot \frac{q+1}{(q_0+1)^3}\cdot \frac{|\mathrm{GU}_n(q_0)|^3}{ |\mathrm{GU}_n(q)|}\\
  &= 36\log_p(q_0)^2 \left(\frac{q_0^3+1}{q_0+1},n\right) \cdot \frac{q_0^3+1}{(q_0+1)^3}\cdot \frac{|\mathrm{GU}_n(q_0)|^3}{ |\mathrm{GU}_n(q_0^3)|}\\
  &>36\cdot\frac{q_0^3+1}{(q_0+1)^3}\cdot \frac{|\mathrm{GU}_n(q_0)|^3}{ |\mathrm{GU}_n(q_0^3)|}.
\end{align*} 
By Lemma~\ref{lm:LUSO} we have 
\[  
\frac{|\mathrm{GU}_n(q_0)|^3}{ |\mathrm{GU}_n(q_0^3)|}>\frac{(1+q_0^{-1})^3(1-q_0^{-2})^3}{(1+q_0^{-3})(1-q_0^{-6})(1+q_0^{-9})}.
\]
Hence, 
\[  
\frac{|H_0|^3|\mathcal{O}_{1}|^2}{|G_0|}>36\cdot \frac{q_0^3+1}{(q_0+1)^3}\cdot \frac{(1+q_0^{-1})^3(1-q_0^{-2})^3}{(1+q_0^{-3})(1-q_0^{-6})(1+q_0^{-9})}=36\cdot \frac{(1-q_0^{-2})^3}{(1-q_0^{-6})(1+q_0^{-9})}:=h(q_0)
\] 
One can easily verify that $h(q_0)$ is increasing in $q_0$, and $h(2)>13$.
Therefore, $|H_0|^3|\mathcal{O}_{1}|^2>|G_0|$ always holds; that is, $H_1$ is always large in $G_1$ when $r=3$. 
\end{proof}



\begin{lemma}\label{lm:psuc6}
Let $G_0=\mathrm{PSU}_n(q)$ and $H_1 \in\mathcal{C}_6$.
Then $H_{1}$ is large in $G_{1}$ if and only if  $H_{0}$ is large in $G_{0}$, and $(G_0,H_0)$ is one of the following:
\[
  (\mathrm{PSU}_{3}(5),3^{2}.\mathrm{Q}_{8} ), \  
 (\mathrm{PSU}_{4}(3),2^{4}.\mathrm{A}_{6} ), \   (\mathrm{PSU}_{4}(7),2^{4}.\mathrm{S}_{6} ). 
\] 
\end{lemma}
\begin{proof}   
By~\cite[Table~4.6.B]{K-Lie},  $H_1$ is of type $r^{1+2m}.\mathrm{Sp}_{2m}(r)$ or   $(4\circ 2^{1+2m}).\mathrm{Sp}_{2m}(2)$.

Assume that $H_1$ is of type  $r^{1+2m}.\mathrm{Sp}_{2m}(r)$, where $n=r^{m},m\geq 1$ and $r$ is an odd prime divisor of $q+1$.
By~\cite[Proposition 4.6.5]{K-Lie}, if $n=3$ and $q\equiv 2,5\pmod{9}$, then $H_0=3^2{:}\mathrm{Q}_8$, otherwise  $H_0=r^{2m}.\mathrm{Sp}_{2m}(r)$.
Then
\[
|H_0|\leq r^{2m}|\mathrm{Sp}_{2m}(r)|\leq  r^{2m+m(2m+1)}=r^{m(2m+3)} .
\]
From Eq.~\eqref{eq:H1} we have 
\[
  \frac{1}{2}q^{n^2-2} <|G_0|\leq  |H_0|^3 \vert \mathcal{O}_{1}\vert^2< r^{3m(2m+3)}(9q^3) .  
\] 
It follows that
\begin{equation}\label{eq:psuc61}
r^{2m}-5=n^2-5 <3m(2m+3)\log_q(r)+\log_{q}(18).
\end{equation}
Since $r \mid (q+1)$, we have $\log_q(r)\leq \log_2(3)<1.6$. Then 
\[
3^{2m}-5<4.8m(2m+3)+\log_{2}(18).
\]
This inequality holds if and only if $m=1$.  
Then  Eq.~\eqref{eq:psuc61} is reduced to 
\[
r^{2}-5<15\log_q(r)+\log_{q}(18).
\]
If $r\geq 5$, then $15\log_q(r)+\log_{q}(18)<15\log_{4}(5)+\log_{4}(18)<20 \leq r^{2}-5$, a contradiction.
Therefore, $r=3$, and hence $n=3$. By~\cite[Table~8.5]{BHRD2013}, $q=p=2\pmod{3}$,   $H_0=3^2{:}\mathrm{Q}_8.c$, where $c=(9,q+1)/3$, and $|\mathcal{O}_{1}|=2(3,q+1)/c$.
Computation shows that Eq.~\eqref{eq:H1} holds if and only if $q=5$, giving the pair $(G_0,H_0)=(\mathrm{PSU}_3(5),3^2{:}\mathrm{Q}_8)$. 

Assume that $H_1$ is of type   $(4\circ 2^{1+2m}).\mathrm{Sp}_{2m}(2)$, where $n=2^m$. 
By~\cite[Proposition 4.6.6]{K-Lie}, $m\geq 2$ and $q=p\equiv 3 \pmod{4}$; moreover, if $n=4$ and $q\equiv 3\pmod{8}$, then $H_0=2^4.\mathrm{A}_6$ and $c=2$, otherwise  $H_0=2^{2m}.\mathrm{Sp}_{2m}(2)$ and $c=(4,q+1)$.

Let $n=4$. Then either  $H_0=2^4.\mathrm{A}_6$, $q=p=3\pmod{8}$ and $|\mathcal{O}_{1}|=4$, or  $H_0=2^4.\mathrm{S}_6$, $q=p=7\pmod{8}$ and $|\mathcal{O}_{1}|=2$.
Computation shows that Eq.~\eqref{eq:H1} holds if and only if $q\in \{3,7 \}$, giving $(\mathrm{PSU}_4(3),2^4.\mathrm{A}_6)$ and $(\mathrm{PSU}_4(3),2^4.\mathrm{S}_6)$.

Let $n>4$. Then $m\geq 3$ and $H_0=2^{2m}.\mathrm{Sp}_{2m}(2)$. Suppose that  $H_1$ is large. Then Eq~\eqref{eq:psuc61} also holds, which implies
\[
2^{2m}-5 <3m(2m+3)\log_q(2)+\log_{q}(18).
\] 
Since  $m\geq 3$, we have $3m(2m+3)\log_q(2)+\log_{q}(18)<81\log_3(2)+\log_{3}(18)<54$, while $2^{2m}-5\geq 59$, a contradiction.
\end{proof}

\begin{lemma}\label{lm:psuc7}
Let $G_0=\mathrm{PSU}_n(q)$ and $H_1 \in \mathcal{C}_7 $.
Then $H_{1}$ is not large in $G_{1}$. 
\end{lemma}

\begin{proof}
Suppose for a contradiction that $H_{1}$ is large in $G_{1}$. 
By~\cite[Table~4.6.B and Proposition 4.7.3]{K-Lie},  $H_1$ is  of type $\mathrm{GU}_{m}(q) \wr \mathrm{S}_t$, where $n=m^t$ with $m\geq 3$ and $t \geq 2$, and $(m,q)\neq (3,2)$, and $|H_0|<|\mathrm{SU}_{m}(q)|^{t}t! $.  
Since $|\mathrm{SU}_{m}(q)|<q^{m^2-1}$, $|G_0|\geq q^{n^2-2}/2$, $\vert \mathcal{O}_{1} \vert^2<9q^3$ and $t!<2^{t\log_2(\frac{t+1}{2})}$, we conclude from Eq.~\eqref{eq:H1} that 
\[
\frac{1}{2}q^{m^{2t}-2} <|G_0|\leq |H_0|^3 \vert \mathcal{O}_{1} \vert^2< (q^{t(m^2-1)}t!)^3\cdot (9q^3)<9q^{3t(m^2-1)+3}\cdot 2^{3t\log_2(\frac{t+1}{2})} .  
\] 
It follows that 
\[
 m^{2t}-3t(m^2-1)-5<\log_2(18) +3t\log_2\left(\frac{t+1}{2}\right).
\]
Since $m\geq 3$, we have
\[
 3^{2t}-24t-5 <\log_2(18) +3t\log_2\left(\frac{t+1}{2}\right). 
\]
However, there is no $t\geq 2$ satisfying the above inequality. 
This contradiction proves the lemma.
\end{proof} 

\subsection{Geometric  subgroups of symplectic groups}\label{sec:proof-psp} 

Throughout this subsection we assume  $G_0=\mathrm{PSp}_n(q)$, where $n \geq 4$ and $(n,q)\neq (4,2)$.
By Table~\ref{tb:orderofsimplegroups},
\[
|\mathrm{Out}(G_0)|=
\begin{cases}
  2\log_p(q),& \text{ if $n=4$ and $q$ is even},\\
   d\log_p(q),& \text{ otherwise},\\  
\end{cases}
\]
where $d=(2,q-1)$. 
We claim that 
\begin{equation}\label{eq:pspout}
\vert \mathcal{O}_{1}\vert^2\leq 4q.
\end{equation} 
If $q=8$, then $d=1$ and  $(\log_p(q))^2=9$, whence $ \vert \mathcal{O}_{1}\vert ^2=9=8q/8<4q$.  
If $q\neq 8$, then $ \log_p(q) ^2\leq q$ by Lemma~\ref{lm:q2f}, and consequently $ \vert \mathcal{O}_{1}\vert ^2\leq 4q$.
Therefore, Eq.~\eqref{eq:pspout}  holds in all cases.

\begin{lemma}[{\cite[Proposition 4.22]{AB2015}}]\label{lm:pspc1}
Let $G_0=\mathrm{PSp}_n(q)$.
If $H_1\in \mathcal{C}_1 \cup \mathcal{C}_8$, then $H_0$ is large in $G_0$. 
\end{lemma}

\begin{lemma}\label{lm:pspc2}
Let $G_0=\mathrm{PSp}_n(q)$ and $H_1 \in\mathcal{C}_2$. Then $H_1$ is large in $G_1$ if and only if one of the following holds:
\begin{enumerate}[\rm (1)]
  \item  $H_1$  is of type $\mathrm{GL}_{n/2}(q).2$.
  \item  $H_1$  is of type $\mathrm{Sp}_{n/t}(q)\wr \mathrm{S}_t$, where  $t\in \{2,3\}$, or $(n,t)=(8,4)$, or $(n,t)= (10,5)$ and $q\in \{3,4,5\}$.  
\end{enumerate} 
In particular, $H_0$ is large in $G_0$ except the case where  $H_1$ is of type $\mathrm{Sp}_{n/t}(q)\wr \mathrm{S}_t$ with $(n,t,q)= (10,5,5)$.
\end{lemma}
\begin{proof}
By~\cite[Table~3.5.C]{K-Lie}, $H_1$ is of type $\mathrm{GL}_{n/2}(q).2 $ or $\mathrm{Sp}_m(q)\wr \mathrm{S}_t$. If $H_1$ is of type $\mathrm{GL}_{n/2}(q).2 $ or $\mathrm{Sp}_{n/t}(q)\wr \mathrm{S}_t$ with $t\in\{2,3\}$, then $H_0$ is large in $G_0$, see ~\cite[Proposition 4.22]{AB2015}. 

Next, we assume that  $H_1$ is of type $\mathrm{Sp}_m(q)\wr \mathrm{S}_t$ with $n=mt$ and $t\geq 4$.
Note that  $(m,q)\neq (2,2)$, see~\cite[Proposition 2.3.6]{BHRD2013}. 
By~\cite[Proposition~4.2.10]{K-Lie}, $|H_0|=d^{-1} |\mathrm{Sp}_m(q)|^t t!$, where $d=(2,q-1)$, and $c=1$ (so $ |\mathcal{O}_{1}| =d\log_p(q)$).  
Notice that $  |G_0|> (2d)^{-1}q^{\frac{n(n+1)}{2}} $ and $|\mathrm{Sp}_m(q)|<q^{\frac{m(m+1)}{2}}$, see Lemmas~\ref{lm:LUSO} and~\ref{lm:LUSO1}, and $\log_p(q)^2<9q/8$.
From Eq.~\eqref{eq:H1} we have   
\[ 
\frac{1}{2d}q^{\frac{tm(tm+1)}{2}} <|G_0|\leq |H_0|^3 \vert \mathcal{O}_{1}\vert^2 < (d^{-3} q^{\frac{3tm(m+1)}{2}} (t!)^3 ) \cdot (\frac{9}{8}q d ^2) . 
\] 
Then
\[
q^{\frac{tm(m(t-3)-2)}{2}-1}<\frac{9}{4}(t!)^3.
\]
Computation shows that  the above inequality holds if and only if  either $(m,t)=(2,4)$ or $(m,t)=(2,5)$ with $q\in \{3,4,5\}$ (note that $(m,q)\neq (2,2)$). 

If $(m,t)=(2,4)$, then $H_0$ is large in $G_0$ by~\cite[Proposition 4.22]{AB2015}. 

Let $(m,t)=(2,5)$ and $q\in \{ 3,4,5\}$. Computation shows that if $q \in \{3,4\}$ then $H_0$ is large in $G_0$; if $q=5$ then $H_1$ is large in $G_1$ while  $H_0$ is not large in $G_0$. 
\end{proof} 
 
\begin{lemma}\label{lm:pspc3}
Suppose $G _0=\mathrm{PSp}_n(q)$ and $H_1 \in\mathcal{C}_3$.
Then $H_{1}$ is large in $G_{1}$ if and only if  $H_0$  is large in $G_0$, and  $H_0$ is of type $\mathrm{GU}_{n/2}(q)$, $\mathrm{Sp}_{n/2}(q^2)$ or $\mathrm{Sp}_{n/3}(q^3)$. 
\end{lemma}

\begin{proof}
By~\cite[Proposition~4.22]{AB2015}, if $H_1$ is of type $\mathrm{GU}_{n/2}(q)$, $\mathrm{Sp}_{n/2}(q^2)$ or $\mathrm{Sp}_{n/3}(q^3)$, then $H_0$ is large in $G_0$.
Hence, it suffices  to prove that $H_{1}$ is not large in $G_{1}$  in the case where $H_0$ is of type $\mathrm{Sp}_{m}(q^r)$ with $r\geq 5$. 
Suppose for a contradiction that this case happen.  
By~\cite[Proposition 4.3.10]{K-Lie}, $|H_0|=d^{-1}r|\mathrm{Sp}_m(q^r)|$ and $c=1$ (so $\mathcal{O}_1=d\log_p^q$).
From Eq.~\eqref{eq:H1} we have  
\[ 
\frac{1}{2d}q^{\frac{rm(rm+1)}{2}} <|G_0|\leq |H_0|^3 \vert \mathcal{O}_{1}\vert^2 <  (d^{-3}r^3 q^{\frac{3rm(m+1)}{2}}  ) \cdot  (\frac{9}{8}q d ^2 ) . 
\]  
Then 
\[
q^{\frac{rm(m(r-3)-2)}{2}-1}<\frac{9}{4}r^3.
\]
Since $m\geq 2$, we derive that
\[
\frac{r(2(r-3)-2)}{2}-1 <\log_2 (\frac{9}{4} )+3\log_2(r).
\]
However, there is no $r\geq 5$ satisfying the above inequality. 
\end{proof}

\begin{lemma}\label{lm:pspc4}
Let $G_0=\mathrm{PSp}_n(q)$ and $H_1\in  \mathcal{C}_4$.
Then $H_{1}$ is not large in $G_{1}$. 
\end{lemma}

\begin{proof}
Suppose for a contradiction that $H_{1}$ is large in $G_{1}$.
By~\cite[Table~3.5.C and Proposition~4.4.11]{K-Lie}, $H_1$ is of type $\mathrm{Sp}_{n_1}(q)\otimes \mathrm{O}_{n_2}^{\epsilon}(q) $, where $q$ is odd, $\epsilon \in \{+,-,\circ \}$, $n_2\geq 3$ and $n={n_1}{n_2}$, and 
\[
|H_0|= |(\mathrm{PSp}_{n_1}(q) \times \mathrm{PO}_{n_2}^{\epsilon}(q)).(2,{n_2})|, \text{ and } c=1.
\] 
Since $q$ is odd, by Proposition~\ref{prop:orders} and Lemma~\ref{lm:LUSO} we have 
\[| \mathrm{PSp}_{n_1}(q) |\leq \frac{1}{2}q^{\frac{1}{2}n_1(n_1+1)} \text{ and } |\mathrm{PO}_{n_2}^{\epsilon}(q)| \leq q^{\frac{1}{2}n_2(n_2-1)}.\] 
Hence, $|H_0|\leq q^{\frac{n_1(n_1-1)+n_2(n_2-1)}{2}}$. 
Since $  |G_0|>  q^{\frac{n(n+1)}{2}}/4 $ and $|\mathcal{O}_1^|2=4\log_p(q)^2\leq  4q$ (as $q$ is odd), we conclude from Eq.~\eqref{eq:H1} that  
\[ 
\frac{1}{4}q^{\frac{{n_1}{n_2}({n_1}{n_2}-1)}{2}} <|G_0|\leq |H_0|^3 \vert \mathcal{O}_{1}\vert^2 < (  q^{\frac{ 3{n_1}({n_1}+1)+3{n_2}({n_2}-1)}{2}} ) \cdot (4q ) . 
\]   
Then $q^{f({n_1},{n_2})}\leq 16$, where 
\[
f({n_1},{n_2}):= \frac{{n_1}{n_2}({n_1}{n_2}-1)-3{n_1}({n_1}-1)-3{n_2}({n_2}-1)}{2}-1.
\]
It is easy to verify that $f({n_1},{n_2})$ is increasing in both $n_1$ and $n_2$.
The only possibility for $(n_1,n_2,q)$ satsifying $q^{f({n_1},{n_2})}\leq 16$ is $(2,3,3)$, where $f(2,3)=2$.
Now $n=6$. However, ~\cite[Table 8.28]{BHRD2013} tells us that $\mathrm{Sp}_{6}(3)$ has no such a maximal subgroup of type $\mathrm{Sp}_{2}(3)\otimes \mathrm{O}_{3} (3) $. 
\end{proof}

\begin{lemma}\label{lm:pspc5}
Let $G_0=\mathrm{PSp}_n(q)$ and $H_1 \in\mathcal{C}_5$. Then $H_{1}$ is large in $G_{1}$ if and only if  $H_1$ is of type $\mathrm{Sp}_n(q^{1/r})$,  where  $r=2$ or $3$. 
In particular,  $H_0$ is large in $G_0$ if and only if  $r=2$. 
\end{lemma}

\begin{proof} 
By~\cite[Table~3.5.C]{K-Lie}, $H_1$ is of type $\mathrm{Sp}_n(q_0 )$, where $q=q_0^r$, $r$ prime.
From~\cite[Proposition~4.22]{AB2015} we see that $H_0$ is large in $G_0$ if and only if $r=2$.
Hence, we assume $r\geq 3$, and we shall prove that $H_1$ is large in $G_1$ only for $r=3$. 
By~\cite[Proposition 4.5.4]{K-Lie}, $H_0=c| \mathrm{PSp}_n(q_0) |$, where $c=(2,q-1,r)$.
Thus, $|H_0|\leq |\mathrm{Sp}_n(q_0)|\leq q_0^{\frac{1}{2}n(n+1)}$.
Since $|G_0|\geq \frac{1}{4}q^{\frac{1}{2}n(n+1)}$ and $|\mathcal{O}_{1}|^2\leq 4q$, by Eq.~\eqref{eq:H1} we have 
\[
  q_0^{\frac{1}{2}rn(n+1)} <|G_0|\leq  |H_0|^3 |\mathcal{O}_{1}|^2< q_0^{\frac{1}{2}3n(n+1)} 4q_0^r .  
\]
It follows that
\[ 
  n(n+1)(r-3)-2r <\log_{q_0}(4).
\] 
If $r\geq 5$, then $n^2(r-3)-2r-3>4^2\cdot 2-10-3=19$, while $\log_{q_0}(4)\leq 2$, a contradiction.
Therefore, $r=3$.
Now,  
\[
  \frac{|H_0|^3|\mathcal{O}_{1}|^2}{|G_0|}=\frac{ c|\mathrm{Out}(G_0)|^2\cdot  |\mathrm{PSp}_n(q_0)|^3}{ |\mathrm{PSp}_n(q_0^3)|} \geq   c\log_p(q)^2 \cdot \frac{|\mathrm{Sp}_n(q_0)|^3}{|\mathrm{Sp}_n(q_0^3)|} .
  \] 
By Lemma~\ref{lm:LUSO}, we have $\frac{1}{2}q_0^{\frac{1}{2}n(n+1)}<|\mathrm{Sp}_n(q_0)|< q_0^{\frac{1}{2}n(n+1)} $.
Hence, ${|\mathrm{Sp}_n(q_0)|^3}/{|\mathrm{Sp}_n(q_0^3)|}\geq \frac{1}{8}$.
Since $3\leq r\leq \log_p(q)$, we have $\log_p(q)^2\geq 9$, and hence $|H_0|^3|\mathcal{O}_{1}|^2>|G_0|$, that is, $H_{1}$ is always large in $G_{1}$. 
\end{proof}

\begin{lemma}\label{lm:pspc6}
Let $G_0=\mathrm{PSp}_n(q)$ and $H_1 \in\mathcal{C}_6$.
Then $H_{1}$ is large in $G_{1}$ if and only if $H_0$ is large in $G_0$, and $(G_0,H_0)$ is one of the following: 
\[
 (\mathrm{PSp}_{4}(3),2^{4}.\mathrm{\Omega}^{-}_{4}(2) ),  
 (\mathrm{PSp}_{4}(5),2^{4}.\mathrm{\Omega}^{-}_{4}(2) ),   (\mathrm{PSU}_{4}(7),2^{4}.\mathrm{O}^{-}_{4}(2)),    (\mathrm{PSp}_{8}(3),2^{6}.\mathrm{\Omega}^{-}_{6}(2)). 
\] 

\end{lemma}
\begin{proof}  
By~\cite[Table~3.5.C and Proposition~4.6.9]{K-Lie}, $H_1$ is of  type $2^{1+2m}.\mathrm{\Omega}_{2m}^{-}(2)$, where $n=2^m$, $q=p\geq 3$, and $c=1$ if $q\equiv \pm 3\pmod{8}$, while $c=2$ if $q\equiv \pm 1\pmod{8}$, and $H_0=2^{2m}.\mathrm{\Omega}_{2m}^{-}(2).c$.
 
If  $m=2$, then computation shows that Eq.~\eqref{eq:H1} holds if and only if $q \in \{3,5,7 \}$, giving pairs $(\mathrm{PSp}_{4}(3),  2^{4}.\mathrm{\Omega}^{-}_{4}(2) )$,   $(\mathrm{PSp}_{4}(5),  2^{4}.\mathrm{\Omega}^{-}_{4}(2)   )$ and $(\mathrm{PSp}_{4}(7),  2^{4}.\mathrm{O}^{-}_{4}(2) )$.

Assume $m\geq 3$. By Lemma~\ref{lm:LUSO}, $|\mathrm{SO}_{2m}^{-}(2)|<  2^{m(2m-1)+1}$. Hence,  $
|H_0|\leq 2^{2m}|\mathrm{SO}_{2m}^{-}(2)|\leq   2^{2m^2+m+2}$. 
Since $  |G_0|>  \frac{1}{4} q^{\frac{2^m(2^m+1)}{2}} $ and $|\mathcal{O}_{1}| \leq 2 $, we conclude from Eq.~\eqref{eq:H1} that  
\[ 
\frac{1}{4}q^{\frac{2^m(2^m+1)}{2}} <|G_0|\leq |H_0|^3 \vert \mathcal{O}_{1}\vert^2 < 2^{6m^2+3m+6} \cdot 4 . 
\]   
It follows that $q^{ 2^{m-1}(2^m+1) }< 2^{6m^2+3m+9}$, implying that
\[
 2^{m-1}(2^m+1) <(6m^2+3m+9)\log_q(2).
\]
If $m\geq 4$, then $2^{m-1}(2^m+1) \geq 136$, while   $(6m^2+3m+9)\log_q(2)  <117$, a contradiction.
Therefore, $m=3$. In the case $m=3$, computation shows that Eq.~\eqref{eq:H1} holds if and only if $q=3$, giving the pair $(\mathrm{PSp}_{8}(3),2^{6}.\mathrm{\Omega}_6^{-}(2))$. 

By~\cite[Proposition 4.22]{AB2015}, $H_0$ is large in $G_0$ for those $4$ pairs $(G_0,H_0)$.
\end{proof}

\begin{lemma}\label{lm:pspc7}
Let $G_0=\mathrm{PSp}_n(q)$ and $H_1 \in \mathcal{C}_7 $.
Then $H_{1}$ is not large in $G_{1}$. 
\end{lemma}

\begin{proof}
Suppose for a contradiction that $H_{1}$ is large in $G_{1}$.
Now, by~\cite[Table 3.5.C and Proposition 4.7.4]{K-Lie}, $H_1$ is  of type $\mathrm{Sp}_{m}(q) \wr \mathrm{S}_t$, where $n=m^t$, and $qt$ is odd, and $(q,m)\neq (3,2)$, and $|H_0|=d^{-1}|\mathrm{Sp}_{m}(q)|^{t}t! $.  
From Eq.~\eqref{eq:H1} we have 
\[ 
\frac{1}{2d}q^{\frac{m^t(m^t+1)}{2}} <|G_0|\leq |H_0|^3 \vert \mathcal{O}_{1}\vert^2 <  (d^{-3} q^{\frac{3tm(m+1)}{2}} (t!)^3  ) \cdot  (\frac{9}{8}q d ^2 ) . 
\] 
Then
\[
q^{\frac{m^t(m^t+1)-3tm(m+1))}{2}-1}<\frac{9}{4}(t!)^3.
\]
A direct computation shows that this inequality fails for every triple  $(q,m,t)$ satsifying  $qt$ odd and $(q,m)\neq (3,2)$.
\end{proof} 


\subsection{Geometric subgroups of orthogonal groups}\label{sec:proof-pso}

Throughout this subsection let \(G_0 = \mathrm{P\Omega}_n^{\epsilon}(q)\) with \(n \ge 7\) and \(\epsilon \in \{+,-,\circ\}\). 
When \(n\) is odd we have \(\mathrm{P\Omega}_n^{\circ}(q) \cong \mathrm{PSp}_{n-1}(q)\) for even \(q\); therefore  we assume that \(q\) is odd whenever \(n\) is odd.

The case \(G_0 = \mathrm{P\Omega}_8^{+}(q)\) is special because the outer automorphism group contains triality automorphisms, which increases \(|\mathcal{O}_1|\) and makes the usual order comparisons more delicate. To avoid this extra complication initially, we first impose the following hypothesis. 

\begin{hypothesis}\label{hy:1}
If \(G_0 = \mathrm{P\Omega}_8^{+}(q)\), we assume that  \(G_1 \) does not  contain any triality automorphism.
\end{hypothesis}

The remaining case where \(G_1\) does contain triality automorphisms will be treated separately in Lemma~\ref{lm:exo8q}. 

Under Hypothesis~\ref{hy:1},  Table~\ref{tb:orderofsimplegroups}  gives
\[
|\mathcal{O}_{1}| \le 
\begin{cases}
2d\log_p(q), & \epsilon \in \{+,-\},\\[4pt]
2\log_p(q),   & \epsilon = \circ,
\end{cases} 
\]
where  $d = (4,\;q^{n/2} - \epsilon)$.
Except for the single exception \((p,q) = (2,8)\), Lemma~\ref{lm:q2f} yields \((\log_p q)^2 \le q\). Using this bound we obtain the estimates
\begin{equation}\label{eq:psoout}
|\mathcal{O}_{1}| ^2\leq \begin{cases}
64q, & \text{if $\epsilon\in\{+,-\}$ and $q$ is odd},\\
9q/2, & \text{if $\epsilon\in\{+,-\}$ and $q $ is even},\\ 
4q, & \text{if $\epsilon=\circ$},\\
\end{cases} 
\end{equation}

\begin{lemma}[{\cite[Proposition 4.23]{AB2015}}]\label{lm:psoc1}
Let $G_0=\mathrm{P\Omega}_n^{\epsilon}(q)$.
If $H_1\in \mathcal{C}_1 $, then $H_0$ is large in $G_0$. 
\end{lemma}

\begin{lemma}\label{lm:psoc2}
Let $G_0=\mathrm{P\Omega}_n^{\epsilon}(q)$, $H_1\in \mathcal{C}_2$ and assume Hypothesis~\ref{hy:1}.
Then $H_{1}$ is large in $G_{1}$ if and only if one of the following holds:
\begin{enumerate}[\rm (1)]
	\item $H_1$ is of type  $\mathrm{GL}_{n / 2}(q). 2$ or $\mathrm{O}_{n / 2}^{\circ}(q)^{2}$. 
	
	\item $H_1$ is of type $\mathrm{O}_{n/2}^{\epsilon_{1}}(q) \wr \mathrm{S}_{2}$, where $\epsilon_{1}\in \{+,-,\circ\}$, and $\epsilon=\epsilon_1^2$ if $n/2 $ is even, while $\epsilon=(-1)^{(q-1)n/4}$ if $n/2 $ is odd.
	\item  The pair $(G_0,H_0)$ is one of the following: 
	\begin{align*}
& (\mathrm{P\Omega}_8^{+}(2),3^4.2^3.\mathrm{S}_4),\ (\mathrm{P\Omega}_8^{+}(3),[2^9].\mathrm{S}_4),\ (\mathrm{P\Omega}_{10}^{-}(2),3^5.2^4.\mathrm{S}_5),\ (\mathrm{P\Omega}_{12}^{-}(2), \mathrm{A}_5^3.2^2.\mathrm{S}_3),\\ 
&(\mathrm{P\Omega}_7(5),2^6.\mathrm{A}_7),\ (\mathrm{P\Omega}_8^{+}(5),2^6.\mathrm{A}_8),\ (\mathrm{P\Omega}_8^{+}(3), 2^{6}.\mathrm{A}_8), \ (\mathrm{P\Omega}_9(3), 2^{8}.\mathrm{A}_9), \ (\mathrm{P\Omega}_{10}^{-}(3), 2^{8}.\mathrm{A}_{10}), \\ 
&  (\mathrm{P\Omega}_{11}(3), 2^{10}.\mathrm{A}_{11}),    \ (\mathrm{P\Omega}_{12}^{-}(3), 2^{10}.\mathrm{A}_{12}),    \ (\mathrm{P\Omega}_{13}(3), 2^{12}.\mathrm{A}_{13}),\ (\mathrm{P\Omega}_{14}^{-}(3),2^{12}.\mathrm{A}_{14}).      
\end{align*}
\end{enumerate} 
In particular, $H_0$ is large in $G_0$  except $(G_0,H_0)=  (\mathrm{P\Omega}_8^{+}(3),[2^9].\mathrm{S}_4)$,  $(\mathrm{P\Omega}_8^{+}(5),2^6.\mathrm{A}_8)$, $(\mathrm{P\Omega}_{14}^{-}(3),2^{12}.\mathrm{A}_{14})$.
\end{lemma}
	
\begin{proof}
By~\cite[Table~2.4]{BHRD2013}, $H_1$ is of one of the following types: $\mathrm{GL}_{n / 2}(q). 2$, $\mathrm{O}_{n / 2}^{\circ}(q)^{2}$, $\mathrm{O}_{1}(p) \wr \mathrm{S}_{n}$,
and  $\mathrm{O}_{m}^{\epsilon_{1}}(q) \wr \mathrm{S}_{t}$.

If $H_1$ is of type  $\mathrm{GL}_{n / 2}(q). 2$ or $\mathrm{O}_{n / 2}^{\circ}(q)^{2}$, then $H_0$ is large in $G_0$, see~\cite[Proposition 4.23]{AB2015}.

Assume that $H_1$ is of type $\mathrm{O}_{m}^{\epsilon_{1}}(q) \wr \mathrm{S}_{t}$ with $m> 1$. Now, $n=mt$. 
By~\cite[Tables~2.4 and 2.5]{BHRD2013}, the following hold:
\begin{itemize}
  \item If $m$ is even then  $\epsilon = \epsilon_1^t$.
  \item If $m$ is odd and $t$ is even, then $\epsilon = (-1)^{(q - 1)n/4}$.
  \item The full preimage of $H_0$ in $\Omega_{n}^{\epsilon}(q)$ is $ \Omega^{\epsilon_1}_m(q)^t.2^{(2,q-1).(t-1)}.\mathrm{S}_t$.
\end{itemize}
Moreover, according to~\cite[Proposition 2.3.6]{BHRD2013}, the maximality of $H_1$ in $G_1$ forces that  $q\geq 5$ when $m=2$, $t=4$ and $\epsilon=+$. 

According to~\cite[Proposition 4.23]{AB2015}, if $t=2$, then $H_0$ is large in $G_0$.
Hence, we assume $t\geq 3$. 


Since $|\mathcal{O}_{1}|^2\leq 64q$ and $|G_0|>\frac{1}{8}q^{\frac{n(n-1)}{2}}$,  we conclude from Eq.~\eqref{eq:H1} that
\begin{equation}\label{eq:psoc2-1}
\frac{1}{8}q^{\frac{mt(mt-1)}{6}-\frac{1}{3}}=\left( \frac{q^{\frac{n(n-1)}{2}}}{8\cdot 64q} \right)^{\frac{1}{3}} <\left( \frac{|G_0|}{|\mathcal{O}|^2} \right)^{\frac{1}{3}} \leq |H_0|. 
\end{equation} 
 
By Lemma~\ref{lm:Omega},   $|\Omega^{-}_2(q)|=\frac{q+1}{(2,q-1)}$, and $|\Omega^{\epsilon_1}_m(q)|<\frac{1}{(2,q-1)}q^{\frac{m(m-1)}{2}}$ for $(m,\epsilon_1)\neq (2,-)$. Then
\[
	\frac{1}{8}q^{\frac{mt(mt-1)}{6}-\frac{1}{3}}< |H_0| \leq \begin{cases}
(q+1)^t\cdot 2^{t-1}\cdot t!, & \text{if } m=2, \\
(q^{\frac{m(m-1)}{2}})^t\cdot 2^{t-1}\cdot t!, & \text{if } m\geq 3.
\end{cases}
\]
The above inequalities hold only when $(q,m,t)$ is one of the following:
\begin{align*}
  &(2,2,4),\ (3,2,4),\ (4,2,4),\ (5,2,4),\ (2,2,5),\ (3,2,5),\ (3,3,3),\ (5,3,3),\\
  &(7,3,3),\ (2,3,4),\ (2,4,3),\ (3,4,3),\ (4,4,3),\ (3,5,3),\ (2,6,3).
\end{align*}

A direct verification shows that $H_{1}$ is large in $G_{1}$ if and only if $(q,m,t,\epsilon,\epsilon_1)=(2,2,4,+,-)$, $(3,2,4,+,-)$, 
  $(2,2,5,-,-)$,  $(2,4,3,-,-)$, giving  $(\mathrm{P\Omega}_8^{+}(2),3^4.2^3.\mathrm{S}_4)$, $(\mathrm{P\Omega}_8^{+}(3),[2^9].\mathrm{S}_4)$, $(\mathrm{P\Omega}_{10}^{-}(2),3^5.2^4.\mathrm{S}_5)$, $(\mathrm{P\Omega}_{12}^{-}(2), \mathrm{A}_5^3.2^2.\mathrm{S}_3)$, respectively.
For these $4$ tuples, $H_0$ is large in $G_0$ except $(G_0,H_0)=(\mathrm{P\Omega}_8^{+}(3),[2^9].\mathrm{S}_4)$, see~\cite[Proposition 4.23]{AB2015}. 

Assume that $H_1$ is of type $\mathrm{O}_{1}(p) \wr \mathrm{S}_{n}$. By~\cite[Table~2.4]{BHRD2013} and~\cite[Proposition~4.2.15]{K-Lie}, the following hold:
\begin{itemize}
  \item $q=p>2$.
  \item $\epsilon=\circ$ iff $n$ is odd; $\epsilon=-$ iff $n \equiv 2\pmod{4}$ and  $q \equiv 3\pmod{4}$;  $\epsilon=+$ otherwise.
  \item $H_0=2^{n-(2,n)-1}.\mathrm{A}_n$ if $q\equiv \pm 3\pmod{8}$. $H_0=2^{n-(2,n)-1}.\mathrm{S}_n$ if $q\equiv \pm 1\pmod{8}$;
  \item $c=2(n,2)$ if $q\equiv \pm 1\pmod{8}$; $c= (n,2)$ if $q\equiv \pm 3\pmod{8}$.
\end{itemize}

Since $|G_0|>\frac{1}{8}q^{\frac{1}{2}n(n-1)}$ and $|\mathcal{O}_{1}|\leq 2d\log_p(q)\leq 8$, we conclude from Eq.~\eqref{eq:H1} that
\[
\frac{1}{8}q^{\frac{1}{2}n(n-1)}<|G_0|\leq |H_0|^3 \vert \mathcal{O}_{1}\vert^2 < (2^{n-1}n!)^3 \cdot 8^2.
\] 
Computation shows that the above inequality holds only if $(n,q)=(7,7)$, or $q=5$ and $7\leq n\leq 9$, or $q=3$ and $7\leq n \leq 15$.

If $(n,q)=(7,7)$, then $|\mathcal{O}_{1}|=1$.
However,  this case is impossible by~\cite[Proposition 4.23]{AB2015}.

If $q=5$ and $n=7$, or $q=3$ and $7\leq n \leq 13$, then $H_0$ is large in $G_0$ by~\cite[Proposition 4.23]{AB2015}, given 
\begin{align*}
&(\mathrm{P\Omega}_7(5),2^6.\mathrm{A}_7),\ (\mathrm{P\Omega}_8^{+}(3), 2^{6}.\mathrm{A}_8), \ (\mathrm{P\Omega}_9(3), 2^{8}.\mathrm{A}_9), \ (\mathrm{P\Omega}_{10}^{-}(3), 2^{8}.\mathrm{A}_{10}), \\ 
&  (\mathrm{P\Omega}_{11}(3), 2^{10}.\mathrm{A}_{11}),    \ (\mathrm{P\Omega}_{12}^{-}(3), 2^{10}.\mathrm{A}_{12}),    \ (\mathrm{P\Omega}_{13}(3), 2^{12}.\mathrm{A}_{13}).      
\end{align*}

For  $(n,q)=(8,5)$, $(9,5)$, $(14,3)$ and $(15,3)$,  we have $(G_0,H_0)=(\mathrm{P\Omega}_8^{+}(5),2^6.\mathrm{A}_8)$, $(\mathrm{P\Omega}_9(5),2^8.\mathrm{A}_9)$, 
$(\mathrm{P\Omega}_{14}^{-}(3),2^{12}.\mathrm{A}_{14})$
 and $(\mathrm{P\Omega}_{15}(3),2^{14}.\mathrm{A}_{15})$, respectively.
Furthermore,  $|\mathcal{O}_{1}|=4,2,4,2$, respectively.
Computation shows that $H_{1}$ is  large in $G_{1}$ if and only if $(G_0,H_0)=(\mathrm{P\Omega}_8^{+}(5),2^6.\mathrm{A}_8)$ or $(\mathrm{P\Omega}_{14}^{-}(3),2^{12}.\mathrm{A}_{14})$.
\end{proof}


%
%
%
%
%
%
%
%

\begin{lemma}\label{lm:psoc3}
Let $G_0=\mathrm{P\Omega}_n^{\epsilon}(q)$, $H_1\in \mathcal{C}_2$ and assume Hypothesis~\ref{hy:1}. Then $H_{1}$ is large in $G_{1}$ if and only if $H_0$ is of type  $\mathrm{O}_{n/2}^{\epsilon}(q^{2})$, $\mathrm{O}_{n/2}^{\circ}(q^{2})$ or $\mathrm{GU}_{n/2}(q)$. In particular, $H_0$ is large in $G_0$.
\end{lemma}

\begin{proof}
By~\cite[Table~2.6]{BHRD2013}, $H_1$ is of one of the following types: $\mathrm{O}_{m}^{\epsilon}(q^{s})$, $\mathrm{O}_{n/2}^{\circ}(q^{2})$ and $\mathrm{GU}_{n/2}(q)$.
If $H_1$ is of type  $\mathrm{O}_{n/2}^{\epsilon}(q^{2})$, $\mathrm{O}_{n/2}^{\circ}(q^{2})$ or $\mathrm{GU}_{n/2}(q)$, then $H_0$ is large in $G_0$ by~\cite[Proposition 4.23]{AB2015}.

Assume that $H_1$ is of type $\mathrm{O}_{m}^{\epsilon}(q^{s})$, where $n=ms$,  $m\geq 3$ and  $s$ is prime.
By~\cite[Table~2.6]{BHRD2013}, 
the preimage of $H_0$ in $\Omega_{n}^{\epsilon}(q)$ is $\Omega_{m}^{\epsilon}(q^s).s$. Hence, $|H_0|\leq  \frac{1}{(2,q-1)}sq^{\frac{1}{2}sm(m-1)}$. Notice that $|\mathcal{O}_{1}|\leq 64q$ when $q$ is odd, and $|\mathcal{O}_{1}|\leq \frac{9}{2}q$ when $q$ is even. 
By Eq.~\eqref{eq:H1} we have 
\[	
\frac{1}{8}q^{\frac{1}{2}n(n-1)}<|G_0|\leq |H_0|^3 \vert \mathcal{O}_{1}\vert^2 <  (\frac{1}{2}sq^{\frac{1}{2}sm(m-1)} )^3 \cdot 64q.
\]
Hence, $q^{\frac{1}{2} ms(ms-1-3(m-1)-1)}<2^{6}s^3q$, which is equivalent to 
\begin{equation}\label{eq:pso-c3-1}
  ms(ms-3m+1)< 2(6\log_q(2)+3\log_q(s)+1). 
\end{equation} 
It is easy to verify that Eq.~\eqref{eq:pso-c3-1} does not hold for $s\geq 5$.
Hence, $s=3$. Computation shows that Eq.~\eqref{eq:pso-c3-1} holds only for $(q,m)=(3, 3)$, $(5,3)$, $(7,3)$, $(2,4)$, $(3,4)$, $(4,4)$, $(3,5)$  or $(2,6)$. 
Further verification shows that Eq.~\eqref{eq:H1} does not hold for any of them.
This completes the proof.
\end{proof}

\begin{lemma}\label{lm:psoc4}
Let $G_0=\mathrm{P\Omega}_n^{\epsilon}(q)$, $H_1\in \mathcal{C}_4$ and assume Hypothesis~\ref{hy:1}.
Then $H_{1}$ is large in $G_{1}$ if and only if  $H_1$ is of type  $\mathrm{Sp}_{2}(q) \otimes \mathrm{Sp}_{n/2}(q)$ and $(n,\epsilon)=(8,+)$ or $(12,+)$. In particular, $H_0$ is large in $G_0$.
\end{lemma}

\begin{proof}
By~\cite[Table 2.7]{BHRD2013}, $H_1$ is of one of the following  types:
\begin{enumerate} [\rm (i)]
  \item $\epsilon=+$, $H_1$ is of type $\mathrm{Sp}_{n_1}(q) \otimes \mathrm{Sp}_{n_2}(q)$, where $n=n_1n_2$ and $n_1 < \sqrt{n}$, and the preimage of $H_0$ in $\Omega_n^{\epsilon}(q)$ is $ (\mathrm{Sp}_{n_1}(q) \circ \mathrm{Sp}_{n_2}(q)).(2,q - 1,n/4)$.
  \item $\epsilon=\pm$, $H_1$ is of type $\mathrm{O}_{n_1}^{\circ}(q) \otimes \mathrm{O}_{n_2}^{\pm}(q) $, where $n=n_1n_2$,  $q$ is odd, $n_1 \geq 3$, $n_2 \geq 4$, and the preimage of $H_0$  in $\Omega_n^{\epsilon}(q)$ is $\mathrm{SO}_{n_1}^{\circ}(q) \times \Omega_{n_2}^{\pm}(q)$.
  \item $\epsilon=\circ$, $H_1$ is of type $ \mathrm{O}_{n_1}^{\circ}(q) \otimes \mathrm{O}_{n_2}^{\circ}(q) $, where $n=n_1n_2$, $3 \leq n_1 < \sqrt{n}$, and the preimage of $H_0$  in $\Omega_n^{\epsilon}(q)$ is $  (\Omega_{n_1}(q) \times \Omega_{n_2}(q)).2$.
  \item $\epsilon=+$, $H_1$ is of type $\mathrm{O}_{n_1}^{\epsilon_1}(q) \otimes \mathrm{O}_{n_2}^{\epsilon_2}(q) $, where $n=n_1n_2$, $n_i \geq 4$ even, $q$ odd, if $\epsilon_1 = \epsilon_2$ then $ n_1 < \sqrt{n}$, and the preimage of $H_0$  in $\Omega_n^{\epsilon}(q)$ is $  (\mathrm{SO}_{n_1}^{\epsilon_1}(q) \circ \mathrm{SO}_{n_2}^{\epsilon_2}(q)).[c]$ (where $c=2$ or $4$, see~\cite[Proposition 4.4.15]{K-Lie}).
\end{enumerate}

Suppose that $H_1$ is as (i). By~\cite[Proposition 4.23]{AB2015}, $H_0$ is large in $G_0$ whenever $n=8$ or $12$. Notice that both $n_1$ and $n_2$ are even. Hence, $n\equiv 0 \pmod{4}$, and we assume $n\geq 16$.
Now $|H_0|\leq |\mathrm{Sp}_{n_1}(q)|\cdot|\mathrm{Sp}_{n_2}(q)|$. Hence, $|H_0|\leq q^{\frac{1}{2}(n_1(n_1+1)+n_2(n_2+1))}$.
By Eq.~\eqref{eq:H1} we have 
\[	
\frac{1}{8}q^{\frac{1}{2}n(n-1)}<|G_0|\leq |H_0|^3 \vert \mathcal{O}_{1}\vert^2 <  (q^{\frac{1}{2}(n_1(n_1+1)+n_2(n_2+1))} )^3 \cdot 64q.
\]
It follows that  
\begin{equation}\label{eq:pso-c4-1}
  n(n-1)-n_1(n_1+1)-n_2(n_2+1)< 2(9\log_q(2)+1)<20. 
\end{equation}  
If $n_1=2$, then $n_2\geq 8$, and Eq.~\eqref{eq:pso-c4-1} becomes $3n_2^2-3n_2-6<20$, which is impossible for any $n_2\geq 8$.
Let $n_1\geq 4$. Now   $n_1,n_2\leq \frac{n}{4}$.
By Eq.~\eqref{eq:pso-c4-1} we have $20> n(n-1)-2\frac{n}{4}(\frac{n}{4}+1)=\frac{1}{8}(7n^2-12n) $.
Since $n\geq 16$, we have $\frac{1}{8}(7n^2-12n)\geq 200 $, a contradiction. Therefore, the case $n\geq 16$ cannot happen.

Suppose that $H_1$ is as  (ii), (iii) or (iv). 
By Lemma~\ref{lm:LUSO}, $|H_0|\leq 4q^{\frac{1}{2}(n_1(n_1-1)+n_2(n_2-1))}$. By Eq.~\eqref{eq:H1}, we obtain that 
\[
  n(n-1)-n_1(n_1-1)-n_2(n_2-1)< 2(15\log_q(2)+1)<32.
\]
 Since $n_1(n_1-1)\cdot n_2(n_2-1)< n_1n_2(n_1n_2-1)$, we have
\[
  n_1n_2(n_1n_2-1)-n_1(n_1-1)-n_2(n_2-1)<(n_1(n_1-1)-1)(n_2(n_2-1)-1).
\]
Hence, 
\[
  (n_1(n_1-1)-1)(n_2(n_2-1)-1)  <32. 
\]
Since $n_1\geq 3$ and $n_2\geq 4$, we see that $(n_1(n_1-1)-1)(n_2(n_2-1)-1)  \geq 55$, a contradiction. 
\end{proof}

\begin{lemma}\label{lm:psoc5}
Let $G_0=\mathrm{P\Omega}_n^{\epsilon}(q)$, $H_1\in \mathcal{C}_5$ and assume Hypothesis~\ref{hy:1}. Then $H_{1}$ is large in $G_{1}$ if and only if $H_0$ is of type $\mathrm{O}_{n}^{\epsilon'}(q_0)$ with $q=q_0^2$ and $\epsilon=\epsilon'^2$, or $\mathrm{O}_{n}^{\epsilon}(q_0)$ with $q=q_0^3$. 
  In particular, $H_0$ is not large in $G_0$ if and only if $H_0$ is of type  $\mathrm{O}_{n}^{\epsilon}(q_0)$ with $q=q_0^3$. 
\end{lemma}

\begin{proof}
From~\cite[Table 4.5.A]{K-Lie}, $H_1$ is of type $\mathrm{O}_{n}^{\epsilon'}(q_0)$ with $q=q_0^r$ and $\epsilon=\epsilon'^r$, where $r$ is a prime.
If $r=2$, then $H_0$ is large in $G_0$ by~\cite[Proposition 4.23]{AB2015}. 

Assume $r\geq 3$. 
Since $r$ is odd and $\epsilon=\epsilon'^r$, we have $\epsilon'=\epsilon$.
According to~\cite[Table 2.8]{BHRD2013}, the preimage of $H_0$ is $\Omega_{n}^{\epsilon}(q_0)$. 
By Eq.~\eqref{eq:H1}, we obtain that
\[
  \frac{1}{8}q_0^{\frac{1}{2}rn(n-1)}<|G_0|\leq |H_0|^3 \vert \mathcal{O}_{1}\vert^2 < q_0^{\frac{3}{2}n(n-1)}  \cdot 64q_0^r.
\]
Hence,
\[
  (r-3)n(n-1)-2r\leq 18\log_{q_0}(2)\leq 18.
\]
If $r\geq 5$, then $(r-3)n(n-1)-2r\geq (r-3)42-2r\geq 74$, a contradiction.
Therefore, $r=3$.

Assume that $n$ is odd. Then $\epsilon=\circ$ and $q$ is odd. Moreover, since $ \Omega_{n}^{\circ} (q) $ has trivial centre, we have $H_0=\Omega_{n}^{\circ} (q_0)$.
Now, $|G_0|=|\mathrm{SO}_{n}^{\circ}(q)|/2 $ and $|H_0|=|\mathrm{SO}_{n}^{\circ}(q_0)|/2$. Thus,
\begin{equation}\label{eq:oc5-1}
 \frac{|H_0|^3}{|G_0|}=\frac{|\Omega_{n}^{\circ}(q_0)|^3}{|\Omega_{n}^{\circ}(q)|} =\frac{|\mathrm{SO}_n^{\circ}(q_0)|^3}{8|\mathrm{SO}_n^{\circ}(q)|}= \frac{\left(|\mathrm{SO}_n^{\circ}(q_0)|/(2,q_0)\right)^3}{8\left(|\mathrm{SO}_n^{\circ}(q)|/(2,q)\right)}. 
\end{equation}

Assume that $n$ is even. Then $\epsilon\in \{+,-\}$, and $\Omega_{n}^{\epsilon} (q)$ has centre of order $z:=\frac{(4,q^m-\epsilon)}{(2,q-1)}$, where $m=n/2$.

If $q$ is odd, then $|H_0|=|\Omega_{n}^{\epsilon}(q_0)|/z $ and $|G_0|=|\Omega_{n}^{\epsilon}(q)|/z$, and hence  

\begin{equation}\label{eq:oc5-2}
 \frac{|H_0|^3}{|G_0|}=\frac{|\Omega_{n}^{\epsilon}(q_0)|^3}{z^2|\Omega_{n}^{\epsilon}(q)|} =\frac{|\mathrm{SO}_n^{\epsilon}(q_0)|^3}{4z^2|\mathrm{SO}_n^{\epsilon}(q)|}= \frac{\left(|\mathrm{SO}_n^{\epsilon}(q_0)|/(2,q_0)\right)^3}{4z^2\left(|\mathrm{SO}_n^{\epsilon}(q)|/(2,q)\right)}. 
\end{equation}

If $q$ is even, then $z=1$, and so 

\begin{equation}\label{eq:oc5-3}
 \frac{|H_0|^3}{|G_0|}=\frac{|\Omega_{n}^{\epsilon}(q_0)|^3}{ |\Omega_{n}^{\epsilon}(q)|} =\frac{|\mathrm{SO}_n^{\epsilon}(q_0)|^3}{4 |\mathrm{SO}_n^{\epsilon}(q)|}= \frac{\left(|\mathrm{SO}_n^{\epsilon}(q_0)|/(2,q_0)\right)^3}{ \left(|\mathrm{SO}_n^{\epsilon}(q)|/(2,q)\right)}. 
\end{equation}

By Lemma~\ref{lm:LUSO}, 
\begin{align*} 
  &(1-q^{-2}-q^{-4}) q^{\frac{1}{2}n(n-1)} <\frac{|\mathrm{SO}_n^{\epsilon}(q)|}{(2,q)} <(1-q^{-2}) q^{\frac{1}{2}n(n-1)}.  
  \end{align*} 
Hence, it follows from Eq.~\eqref{eq:oc5-1}--Eq.~\eqref{eq:oc5-3} that 
\[
kf(q_0)  \leq \frac{|H_0|^3}{|G_0|} \leq  kg(q_0),
\]
where $k \in \{4,8,16\}$ if $q$ is odd, or $1$ if $q$ is even, and 
\[
f(q_0)=\frac{(1-q_0^{-2}-q_0^{-4})^3 }{1-q_0^{-6} },\ g(q_0)=\frac{(1-q_0^{-2} )^3 }{1-q_0^{-6}-q_0^{-12}}.
\]  
Computation shows that $0.33<f(q_0)<g(q_0)<1$.
By~\cite[Tables 3.5.E and 3.5.F]{K-Lie}, $|\mathcal{O}_{1}|=|\mathrm{Out}(G_0)|=2d\log_p(q)$.
Hence, $|\mathcal{O}_{1}|\geq 2r\geq 6$ if $q$ is even, or $|\mathcal{O}_{1}|\geq 4r\geq 12$ if $q$ is odd.
It follows that  $\frac{|H_0|^3|\mathcal{O}_{1}|^2}{|G_0|}\geq k\cdot 0.33\cdot 6^2 >1 $, implying that $H_{1}$ is large in $G_{1}$. 
\end{proof}

\begin{lemma}\label{lm:psoc6}
Let $G_0=\mathrm{P\Omega}_n^{\epsilon}(q)$, $H_1\in \mathcal{C}_6$ and assume Hypothesis~\ref{hy:1}.
Then $H_{1}=H_0$ and $G_{1}=G_0$, and $H_{1}$ is large in $G_{1}$ if and only if  $(G_0,H_0)=(\mathrm{P\Omega}_8^{+}(3),2^6.\Omega_{6}^+(2))$.  
\end{lemma}

\begin{proof}
By~\cite[Tables 3.5.E and 3.5.G]{K-Lie}, $\mathcal{O}_{1}=1$, that is, $H_{1}=H_0$ and $G_{1}=G_0$.
Hence, it follows from~\cite[Proposition 4.23]{AB2015} that $H_{1}$ is large in $G_{1}$ if and only if $(G_0,H_0)=(\mathrm{P\Omega}_8^{+}(3),2^6.\Omega_{6}^+(2))$. 
\end{proof}

\begin{lemma}\label{lm:psoc7}
Let $G_0=\mathrm{P\Omega}_n^{\epsilon}(q)$, $H_1\in \mathcal{C}_6$ and assume Hypothesis~\ref{hy:1}.
Then $H_{1}$ is not large in $G_{1}$.  
\end{lemma}

\begin{proof}
By~\cite[Table~2.10]{BHRD2013}, there are three cases:
\begin{enumerate} [\rm (i)]
\item $\epsilon=+$,  $H_1$ is of type  $\mathrm{Sp}_{m}(q)\wr\mathrm{S}_{t}$, where $n=m^t$, $qt$ is even and  $(m,q)\notin\{(2,2),(2,3)\}$; 
\item $\epsilon=\circ$, $H_1$ is of type  $\mathrm{O}_{m}(q)\wr\mathrm{S}_{t}$, where $n=m^t$, $m\geq3$ and $(m,q)\neq(3,3)$;
\item  $\epsilon=+$,  $H_1$ is of type  $\mathrm{O}_{m}^{\epsilon'}(q)\wr\mathrm{S}_{t}$, where $n=m^t$,  $q $ is odd, and $m\geq 4$ if $\epsilon' = +$, and $m\geq 6$ if $ \epsilon' = -$.
\end{enumerate}

Suppose that (i) happens. From~\cite[Table 8.50]{BHRD2013}, we see that, if $n=8$, then $G_{1}$ has no such maximal subgroups $H_{1}$ of type $\mathrm{Sp}_{2}(q)\wr\mathrm{S}_{3}$. Hence, we assmue that $(m,t)\neq (2,3)$. 

By~\cite[Proposition 4.7.5]{K-Lie}, $|H_0|\leq |\mathrm{PSp}_m(q)|^t 2^{t-1}\mathrm{S}_t$.
Notice that $|\mathrm{PSp}_m(q)|\leq q^{\frac{1}{2}m(m+1)}$, $|G_0|\geq \frac{1}{8}q^{\frac{1}{2}n(n-1)}$, $t!\leq 2^{t\log_2(\frac{t+1}{2})}$ and $|\mathcal{O}_{1}|^2\leq 64q$.
It follows from Eq.~\eqref{eq:H1} that
\[
\frac{1}{8}q^{\frac{1}{2}n(n-1)}<|G_0|\leq |H_0|^3|\mathcal{O}_{1}|^2\leq (q^{\frac{1}{2}mt(m+1)}2^{t-1}2^{t\log_2(\frac{t+1}{2})})^3\cdot 64q.
\]
Then we derive  that
\begin{equation}\label{eq:psoc7-1}
m^t(m^t-1)-3mt(m+1)-1<\log_q(2)( 6+3t +3t\log_2(\frac{t+1}{2})).   
\end{equation} 
It is easy to see that $m^t(m^t-1)-3mt(m+1)-1$ is increasing on $m$.
Therefore, 
\[
 2^t(2^t-1)-18t-1 \leq m^t(m^t-1)-3mt(m+1)-1< 6+3t +3t\log_2(\frac{t+1}{2}) .   
\]
It is easy to verify that there is no $t\geq 4$ satisfying the above inequality. Hence, case (i) cannot happen.

Suppose that (ii) happens. 
By~\cite[Proposition 4.7.8]{K-Lie}, $|H_0|=|\Omega_{m} (q)^{t}.2^{t - 1}.\mathrm{S}_{t}|$.
Notice that $|\Omega_{m} (q)|<q^{\frac{1}{2}m(m-1)} $.
By Eq.~\eqref{eq:H1}, we  have
\begin{equation}\label{eq:psoc7-2}
  m^t(m^t-1)-3mt(m-1)-1<\log_q(2)(6+ 3t+3t\log_2(\frac{t+1}{2})).   
\end{equation} 
Now both $m $ and $q$ are odd and greater than $2$. Then we have
\[
3^t(3^t-1)-18t-1\leq m^t(m^t-1)-3mt(m-1)-1< 6+ 3t+3t\log_2(\frac{t+1}{2}).
\]
It is easy to verify that there is no $t\geq 3$ satisfying the above inequality. 

Suppose that (iii) happens. By~\cite[Proposition 4.7.5]{K-Lie}, $|H_0|\leq |\mathrm{SO}_m^{\epsilon'}(q)|^t 2^{t-1}\mathrm{S}_t$.
Since $|\mathrm{SO}_m^{\epsilon'}(q)|\leq q^{\frac{1}{2}m(m-1)}$ by Lemma~\ref{lm:LUSO} (now $q$ is odd), we see that Eq.~\eqref{eq:psoc7-2} also holds here. Notice that $m\geq 4$.
Hence, this case also cannot happen. 
\end{proof}

\begin{lemma}\label{lm:exo8q}
Let $G_0=\mathrm{P\Omega}_8^{+}(q)$ and $H_1\in\mathcal{C}_1 \cup \cdots \cup \mathcal{C}_8$, and suppose that $G_{1}$ contains triality automorphisms.
Then $H_{1}$ is large in $G_{1}$ if and only if either $H_0$ is large in $G_0$, or $H_1$ is a $\mathcal{C}_5$-subgroup of type $\mathrm{O}_{8}^{+}(q_0)$ with $q=q_0^3$. 
\end{lemma} 
\begin{proof} 
By~\cite[Table~8.50]{BHRD2013} the candidates for $H_1$ are 
\begin{enumerate}[\rm (i)]
 
\item  $\mathcal{C}_1$-subgroups of type $P_2$ (parabolic subgroups); 
\item  $\mathcal{C}_1$-subgroups of type $\mathrm{O}_{2}^{+}(q)\perp \mathrm{O}_{6}^{+}(q)$ and $\mathrm{O}_{2}^{-}(q)\perp \mathrm{O}_{6}^{-}(q)$;
\item $\mathcal{C}_2$-subgroups of type $\mathrm{O}_{4}^{+}(q)\wr \mathrm{S}_2$ with $q\geq 3$;

\item $\mathcal{C}_2$-subgroups of type $\mathrm{O}_{2}^{-}(q)\wr \mathrm{S}_4$;
\item $\mathcal{C}_2$-subgroups of type $\mathrm{O}_{2}^{+}(q)\wr \mathrm{S}_4$ with $q\geq 5$;

\item  $\mathcal{C}_5$-subgroups of type $\mathrm{O}_{8}^{+}(q_0)$ with $q=q_0^r$, $r$ prime. 
\end{enumerate}

Notice that $|G_0|>q^{28}/8$  and  $|\mathcal{O}_{1} |^2\leq (24\log_p(q))^2\leq 24^2\cdot 9q/8=648q$.
 
If $H_1$ satisfies case (i), (ii) or (iii), then $H_0$ is large in $G_0$ by~\cite[Proposition 4.23]{AB2015}.

Suppose that $H_1$ is as (iv). By~\cite[Table 8.50]{BHRD2013}, the full preimage of $H_0$ in $\Omega_{8}^{+}(q)$ is $\Omega_{2}^{-}(q)^4.(2d)^3.\mathrm{S}_4$
Thus, $|H_0|=(\frac{q+1}{2})^4\cdot   2^3\cdot (2,q-1)^2\cdot (4!)<48\cdot (\frac{3}{2})^4 \cdot q^4< 243 q^4$.
Since $|\mathcal{O}_{1} |^2 \leq 648q$, from Eq.~\eqref{eq:H1} we derive that $243^3\cdot 648\cdot  8>q^{23}$.
This inequality holds only for $q=2$.
From~\cite[Proposition 4.23]{AB2015} that $H_0$ is large in $G_0$ when $q=2$.

Suppose that $H_1$ is as (v). Now the full preimage of $H_0$ in $\Omega_{8}^{+}(q)$ is $\Omega_{2}^{+}(q)^4.(2d)^3.\mathrm{S}_4$, where $d=(2,q-1)$.
Since $|\Omega_{2}^{+}(q)|=(q-1)/d<(q+1)/d=|\Omega_{2}^{-}(q)|$, by the argument for case (iv), we  only need to consider the case $q=2$. This is a contradiction to the condition $q\geq 5$ in case (v). Therefore, in this case $H$ is not large.

Suppose that $H_1$ is as (vi).
If $r=2$, then $H_0$ is large in $G_0$.
Assume $r\geq 3$.
By the argument in Lemma~\ref{lm:psoc5}, we see that if $r\geq 5$ then $H_1$ is not large in $G_1$, and if $r=3$, then
\(
0.33k \leq  {|H_0|^3}/{|G_0|} \leq  k
\),
where $k \in \{4,8,16\}$ if $q$ is odd, or $1$ if $q$ is even.
Let $r=3$.
Now, $|\mathcal{O}_1|=24\log_p(q)$ if $q$ is odd, or $|\mathcal{O}_1|=6\log_p(q)$ if $q$ is even.
Since $0.33\cdot 16^{-1}\cdot 24^2>1$ and $0.33\cdot 6^2>1$, we conclude that $H_1$ is large in $G_1$.
\end{proof}

\subsection{Subgroups in the collection $\mathcal{S}$}\label{sec:proof-as} 
 
We now turn to the case where $H_1 \in \mathcal{S}$. 

Let $V$ be the natural module of $G_0$ and set $n=\dim(V)$. 
Since $H_1 \in \mathcal{S}$, the socle $\mathrm{soc}(H_1) $ is a nonabelian simple group, and the projective action of $\mathrm{soc}(H_1)$ on $V$ corresponds to an absolutely irreducible representation of the covering group of $\mathrm{soc}(H_1)$ in $\mathrm{GL}(V)$.   

Following~\cite{AB2015}, we define a special subcollection $\mathcal{A}$ of $\mathcal{S}$ consisting of those subgroups $H_1$ for which $ \mathrm{soc}(H_1)=\mathrm{A}_{d}$, and  $V$ is the fully deleted permutation module for $\mathrm{A}_{d}$ over $\mathbb{F}_p$, and the group $G_0$ belongs to one of the families listed in Table~\ref{tab:A0}.

\begin{table}[ht]
\centering
\caption{The groups $G_0$ for the subcollection $\mathcal{A}$ of $\mathcal{S}$.}
\label{tab:A0}
\begin{tabular}{ c c c  }
\hline
$d$ & $p$ & $G_0$  \\ \hline
$d\equiv2\pmod{4}$ & 2 & $\mathrm{Sp}_{d-2}(2)$   \\  
$d\equiv0\pmod{8}$ & 2 & $\mathrm{P\Omega}_{d-2}^{+}(2)$\\
$d\equiv 4\pmod{8}$ & 2 & $\mathrm{P\Omega}_{d-2}^{-}(2)$ \\
$d\equiv\pm 1\pmod{8}$ & 2 & $\mathrm{P\Omega}_{d-1}^{+}(2)$\\
$d\equiv \pm 3\pmod{8}$ & 2 & $\mathrm{P\Omega}_{d-1}^{-}(2)$  \\
 $(d,p)=1$ & odd & $\mathrm{P\Omega}_{d-1}^{\epsilon}(p) $ \\
 $(d,p)\neq 1$ & odd & $\mathrm{P\Omega}_{d-2}^{\epsilon}(p) $ \\
 \hline
\end{tabular}
\end{table}

\begin{lemma}
Suppose $H_0\in\mathcal{A}$.
Then $H_{1}$ is large in $G_{1}$ if and only if $H_0$ is large in $G_0$ and $(G_0,H_0)$ lies in Table~\ref{tab:S}.
\end{lemma}

\begin{proof} 
Assume that \(H_{1}\) is large in \(G_{1}\). Then \(|H_{1}|^3 \ge |G| \ge |G_0|\). 
Using Lemma~\ref{lm:LUSO1} we have the lower bound  $|G_0| >\frac{1}{8}q^{\frac{1}{2}n(n-1)}\geq  2^{\frac{1}{2}n(n-1)-3}$.
By Lemma~\ref{lm:n!}, $|H_{1}|\leq|\mathrm{S}_d| =d!<2^{d\log_2(\frac{d+1}{2}) }$.
 Combining these estimates yields 
\[
2^{\frac{1}{2}n(n-1)-3}<|G_0|\leq |H_{1}|^3 \leq 2^{3d\log_2(\frac{d+1}{2})}.
\] 
Since $n\geq d-2$ by Table~\ref{tab:A0}, we have 
\[
\frac{1}{2}(d-2)(d-3)-3<3d\log_2(\frac{d+1}{2}).
\] 
A direct verification shows that this inequality holds only for $d\leq 28$. 
For each $d\in \{5,6,\ldots,28\}$, we compute the exact value of $(d!)^3$ and $|G_0|$ with $(d,G_0)$ satisfying Table~\ref{tab:A0}. Our calculations reveal that $(d!)^3\geq |G_0|$ holds only when $d\leq 24$ in the case $q=2$, and when $d\leq 12$ for $q\geq 3$. 
The precise pairs \((G_0,H_0)\) satisfying these bounds are listed in \cite[Section~8.2]{BHRD2013} (for \(d \le 12\)) and in \cite[Table~9]{AB2015} (for \(12 < d \le 24\)).
A case-by-case verification confirms that for every such pair we indeed have \(|H_0|^3 \ge |G_0|\); consequently all these pairs belong to \cite[Table 9]{AB2015}, which we have reproduced in our own Table~\ref{tab:S}. This completes the proof.
\end{proof}

\begin{table}[h]
\centering
\caption{The pairs $(G_0,H_0)$ with $H_1 \in \mathcal{S} $ such that $H_1$ is large in $G_1$.}
\label{tab:S}
\begin{tabular}{lllllll} 
\hline 
$G_0$ & $H_0$ & Conditions  & $H_1 \in \mathcal{A}$? & $|H_0|^3\geq |G_0|$? \\   \hline 
$\mathrm{P\Omega}_{22}^{+}(q)$ & $\mathrm{A}_{24}$ & $q=2$& Yes& Yes\\
$\mathrm{PSp}_{20}(q)$ & $\mathrm{S}_{22}$ & $q=2$& Yes& Yes\\
$\mathrm{P\Omega}_{20}^{-}(q)$ & $\mathrm{A}_{21}$ & $q=2$& Yes& Yes\\
$\mathrm{P\Omega}_{18}^{-}(q)$ & $\mathrm{A}_{20}$ & $q=2$& Yes& Yes\\
$\mathrm{PSp}_{16}(q)$ & $\mathrm{S}_{18}$ & $q=2$& Yes& Yes\\
$\mathrm{P\Omega}_{16}^{+}(q)$ & $\mathrm{A}_{17}$ & $q=2$& Yes& Yes\\
$\mathrm{P\Omega}_{14}^{+}(q)$ & $\mathrm{A}_{16}$ & $q=2$& Yes& Yes\\
$\mathrm{PSp}_{12}(q)$ & $\mathrm{S}_{14}$ & $q=2$& Yes& Yes\\
$\mathrm{P\Omega}_{12}^{-}(q)$ & $\mathrm{A}_{13}$ & $q=2$& Yes& Yes\\
$\mathrm{P\Omega}_{10}^{+}(q)$ & $\mathrm{A}_{12}$ & $q=3$& Yes& Yes\\
$\mathrm{P\Omega}_{10}^{-}(q)$ & $\mathrm{A}_{12}$ & $q=2$& Yes& Yes\\
$\mathrm{P\Omega}_{10}^{-}(q)$ & $\mathrm{M}_{12}$ & $q=2$& No& Yes\\
$\mathrm{PSU}_9(q)$ & $ \mathrm{J}_3$ & $ q=2$ & No& Yes  \\ 
$\mathrm{P\Omega}_9(q)$ & $ \mathrm{A}_{10}$ & $ q=3$ & Yes& Yes  \\ 
$\mathrm{PSp}_8(q)$ & $ \mathrm{S}_{10}$ & $ q=2$  & Yes& Yes \\ 
$\mathrm{P\Omega}_{8}^{+}(q)$ & $\mathrm{P\Omega}_8^{+}(2)$ & $q \in \{3, 5, 7\}$ & No& Yes \\    
 & $\mathrm{A}_{10}$ & $q=5$ & No& No\\   
  & $\mathrm{P\Omega}_{7}(q)$ &  & No& Yes \\ 
    & $\mathrm{P\Omega}_{8}^{-}(q_0) $ & $q=q_0^2$ & No& Yes  \\ 
  &  ${}^3\mathrm{D}_4 (q_0) $ & $q=q_0^3$ and $q$ is odd &   No& Yes \\ 
   &  ${}^3\mathrm{D}_4 (q_0) $ & $q=q_0^3$ and $q$ is even &  No& No    \\     
$\mathrm{P\Omega}_7 (q)$ & $\mathrm{S}_9$ & $q =3$ & Yes& Yes\\  
 & $\mathrm{PSp}_6(2)$ & $q \in \{3, 5, 7\}$ & No& Yes\\  
   & $\mathrm{G}_{2}(q) $ &  $q$ is odd  & No& Yes\\  
$\mathrm{PSU}_{6}(q)$ & $\mathrm{M}_{22}$ & $q = 2$& No& Yes \\  
& $\mathrm{PSU}_{4}(3).2$ & $q = 2$ & No& Yes\\

$\mathrm{PSp}_{6}(q)$ & $\mathrm{PSU}_{3}(3).2$ & $q = 2$ & No& Yes\\ 
 & $\mathrm{G}_{2}(q) $ & $q$ is even & No& Yes\\  
& $\mathrm{J}_{2} $ & $q = 4$ & No& Yes\\

$\mathrm{PSL}_{5}(q)$ & $\mathrm{M}_{11}$ & $q = 3$ & No& Yes\\  
$\mathrm{PSU}_{5}(q)$ & $\mathrm{PSL}_{2}(11)$ & $q = 2$ & No& Yes\\  

$\mathrm{PSL}_4(q)$ & $\mathrm{A}_7$ & $q=2$  & No& Yes\\
   & $\mathrm{PSU}_4(2)$ & $q=7$  & No& Yes \\
$\mathrm{PSU}_4(q)$ & $\mathrm{A}_7$ & $q\in \{3,5\}$  & Yes& Yes\\
   & $\mathrm{PSU}_4(2)$ & $q=5$ & No& Yes \\
 & $\mathrm{PSL}_3(4)$ & $q=3$ & No& Yes \\ 
 
 $\mathrm{PSp}_{4}(q)$ & $\mathrm{A}_{5}$ & $q = 2$ &   No& Yes \\    
 & $\mathrm{A}_{6}$ & $q = 5$ & Yes& Yes\\
    & $\mathrm{A}_{7}$ & $q = 7$ & Yes& Yes\\ 
  & ${}^2\mathrm{B}_2(q)$ & $q$ is even  & No& Yes\\ 
$\mathrm{PSL}_3(q)$ & $\mathrm{A}_6$ & $q=4$ & No& Yes\\ 
$\mathrm{PSU}_3(q)$ & $\mathrm{A}_6$ & $q=11$ & No& No \\ 
  & $\mathrm{A}_7$ & $q=5$ & No& Yes\\ 
& $\mathrm{M}_{10}$ & $q=5$ & No& Yes\\ 
& $\mathrm{PSL}_{2}(7)$ & $ q\in\{3,5\}$ & No& Yes  \\  
$\mathrm{PSL}_2(q)$ & $\mathrm{A}_5$ & $q\in \{11, 19, 29, 31, 41, 59, 61, 71 \}$ & No& Yes\\ 
& $\mathrm{A}_5$ & $q\in \{9,49 \}$ & Yes& Yes\\ \hline
\end{tabular}
\end{table}

\begin{lemma}
Suppose $H_0\in\mathcal{S}\setminus \mathcal{A}$.
Then $H_{1}$ is large in $G_1$ if and only if  $(G_0,H_0)$ lies in Table~\ref{tab:S}. 
\end{lemma}

\begin{proof} 
By Liebeck~\cite[Theorem~4.2]{L1985}, either $(G_0,\mathrm{soc}(H_1))$ satisfies~\cite[Table 4]{L1985} or $|H_{1}|<q^{2n+4}$. 
Moreover, \cite[Theorem]{L1985} gives the universal bound  $|H_{1}|<q^{3n}$.
Combining this bound  with the inequality $|H_{1}|^3\geq | G_0|> \frac{1}{8}q^{\frac{1}{2}n(n-1)}$ we deduce $q^{9n}\geq \frac{1}{8}q^{\frac{1}{2}n(n-1)}$, and hence  
\[
  18n\geq n(n-1)-6\log_q(2)\geq  n(n-1)-6.
\]
It holds only for $n<20$ (also see the proof of~\cite[Corollary~4.31]{AB2015}).

Assume that $(G_0,\mathrm{soc}(H_{1}))$ satisfies~\cite[Table 4]{L1985}. Since $n<20$, we have only a few possible cases for $(G_0,\mathrm{soc}(H_{1}))$: 
\[
\text{$(\mathrm{PSL}_{\frac{1}{2}d(d-1) }(q),\mathrm{PSL}_{d}(q))$, $(\mathrm{P\Omega}_8^{+}(q),\mathrm{P\Omega}_7(q))$, $(\mathrm{P\Omega}_{16}^{+}(q),\mathrm{P\Omega}_9(q))$, $(\mathrm{PSL}_{11}(2),\mathrm{M}_{24})$}.
\] 
A direct verification shows that Eq.~\eqref{eq:H1} holds only when 
\[(G_0,\mathrm{soc}(H_{1})) \in \{(\mathrm{PSL}_{6}(q),\mathrm{PSL}_{4}(q)),(\mathrm{P\Omega}_8^{+}(q),\mathrm{P\Omega}_7(q))\}.\]
Since  $\mathrm{PSL}_{4}(q)\cong \mathrm{P\Omega}_6^{+}(q)$ is viewed as a $\mathcal{C}_8$-subgroup of $\mathrm{PSL}_{6}(q)$ (see~\cite[Tables 8.24 and 8.25]{BHRD2013}), the pair $(\mathrm{PSL}_{6}(q),\mathrm{PSL}_{4}(q))$ is not listed in Table~\ref{tab:S}.
One checks easily that $\mathrm{P\Omega}_7(q)$ is large in $\mathrm{P\Omega}_8^{+}(q)$, and so this pair $(\mathrm{P\Omega}_8^{+}(q), \mathrm{P\Omega}_7(q))$ is listed in Table~\ref{tab:S}.

Assume $|H_{1}|<q^{2n+4}$. Together with the inequality $|H_{1}|^3\geq | G_0|> \frac{1}{8}q^{\frac{1}{2}n(n-1)}$,  we derive that
\[
12n+24>n(n-1)-6.
\]  
This holds only when $n\leq 14$. For all classical groups with $n\leq 14$ the complete lists of maximal subgroups are available in~\cite{BHRD2013} and~\cite{Schr2015}. 
Examining every possible pair $(G_{0},H_{0})$ with $H_0\in \mathcal{S}\setminus \mathcal{A}$ and checking the condition $|H_1|^3\geq |G_1|$ leaves exactly those pairs recorded in Table~\ref{tab:S}. 
\end{proof} 
 
\begin{remark}\label{remark:3d4q}
Among all pairs \((G_0, H_0)\) with \(H_1 \in \mathcal{S} \setminus \mathcal{A}\) and \(n \leq 14\), all cases are straightforward except for
\(
(G_0, H_0) = \bigl(\mathrm{P\Omega}_8^{+}(q), {}^3\mathrm{D}_4(q_0)\bigr)\) with $q=q_0^3$, which needs a bit more care. When \(q\) is odd,  \(H_0\) is already large in \(G_0\) (see \cite[Table 7]{AB2015}). For even \(q\), using the technique described in Subsection~\ref{sec:techniques}, one can obtain the estimate
\(
1 < \left( {|G_0|}/{|H_0|^3}\right)^{\!1/2} < 2,
\)
which implies that \(H_1\) is indeed large in \(G_1\) (note that $|\mathcal{O}_1|=3\log_p(q)$, see~\cite[Table 8.50]{BHRD2013}).
\end{remark}

\subsection{Subgroups in the collection $\mathcal{N}$} \label{sec:proof-N}

In this subsection, we consider those subgroups $H_1 \in  \mathcal{N}$. 

As noted in~\cite[Section 2.6]{BG2016}, $G_0=\mathrm{PSp}_4(q)$ with $q=2^e\geq 4$ (note that we assume $(n,q)\neq (4,2)$ when $G_0$ is a symplectic group in Theorem~\ref{th:psp}), or $G_0=\mathrm{P\Omega}_8^{+}(q)$.
Moreover, $G_1$ contains graph-field or triality automorphisms, respectively, and the pairs $(G_0,H_0)$ are listed in~\cite[Table~5.9.1]{BG2016}.

\begin{lemma}
Suppose $H_1\in\mathcal{N} $ and $G_0=\mathrm{PSp}_4(q)$ with $q=2^e\geq 4$.
Then $H_{1}$ is large in $G_1$ if and only if either $H_0=[q^4]{:}(q-1)^2$, or 
\[
(G_0,H_0)\in \{ (\mathrm{PSp}_4(8), 7^2{:}\mathrm{D}_8), (\mathrm{PSp}_4(4), 5^2{:}\mathrm{D}_8), (\mathrm{PSp}_4(8), 9^2{:}\mathrm{D}_8), (\mathrm{PSp}_4(4), 17{:}4)\}. 
\] 
In particular, $H_0$ is not large in $G_0$ if and only if  $(G_0,H_0)=(\mathrm{PSp}_4(8), 7^2{:}\mathrm{D}_8)$, $(\mathrm{PSp}_4(8), 9^2{:}\mathrm{D}_8)$ or $ (\mathrm{PSp}_4(4), 17{:}4)$.
\end{lemma}

\begin{proof} 
According to~\cite[Table~5.9.1]{BG2016}, the possible types of $H_0$ are $\mathrm{O}_2^{\epsilon}(q)\wr \mathrm{S}_2$,  $\mathrm{O}_2^{-}(q^2).2$ and $P_{1,2}=[q^4].\mathrm{GL}_1(q)$.
By~\cite[Table~8.14]{BHRD2013}, these correspond to the concrete groups $(q-1)^2{:}\mathrm{D}_8$, $(q+1)^2{:}\mathrm{D}_8$, $(q^2+1){:}4$  and $[q^4]{:}(q-1)^2$, and in every case we have $|\mathcal{O}_1|= |\mathrm{Out}(G_0)|=2\log_2(q) $. 
Note that $|G_0|=q^4(q^4-1)(q^2-1)$. 
If $H_0$ is $[q^4]{:}(q-1)^2$, then clearly $H_0$ is large in $G_0$.
 
Assume $H_0=(q-1)^2{:}\mathrm{D}_8$. From \cite[Table~8.14]{BHRD2013} we see that this group occurs only for $q>4$.  
Using Lemma~\ref{lm:LUSO} and Lemma~\ref{lm:q2f}, we conclude from Eq.~\eqref{eq:H1} that
\[
(1-q^{-2}-q^{-4})q^{10}<|G_0|\leq | H_0|^3 \cdot |\mathcal{O}_1|^2=(8(q-1)^2)^3\cdot (2\log_2(q))^2<2^{11}q^7.
\]
This forces $q=8$.
For $q=8$, a direct verification confirms that $H_1$ is indeed large in $G_1$.

The cases $H_0=(q+1)^2{:}\mathrm{D}_8$ and $H_0=(q^2+1){:}4$ are handled by completely analogous computation. 
\end{proof} 

\begin{lemma}
Suppose $H_1\in\mathcal{N} $ and $G_0=\mathrm{P\Omega}_8^{+}(q)$.
Then $H_{1}$ is large in $G_1$ if and only if $H_0$ is of type $\mathrm{GL}_1(q)\times \mathrm{GL}_3(q)$, $\mathrm{GU}_1(q)\times \mathrm{GU}_3(q)$, $ \mathrm{G}_2(q)$, $P_{1,3,4}$ (a non-maximal parabolic subgroup), or $(G_0,H_0)\in\{(\mathrm{P\Omega}_8^{+}(2), \mathrm{D}_{10}^2.2.\mathrm{S}_2),(\mathrm{P\Omega}_8^{+}(3),[2^9].\mathrm{SL}_3(2))\}$.  
In particular, $H_0 $ is not large in $G_0$ if and only if $(G_0,H_0)=(\mathrm{P\Omega}_8^{+}(2), \mathrm{D}_{10}^2.2.\mathrm{S}_2)$ or $(\mathrm{P\Omega}_8^{+}(3),\mathrm{PSL}_3(3).2)$. 
\end{lemma}

\begin{proof}
According to \cite[Table~5.9.1]{BG2016}, the possible types for \(H_0\) are:
\[
\mathrm{GL}_1^{\epsilon}(q) \times \mathrm{GL}_3^{\epsilon}(q),\;
\mathrm{O}_2^{-}(q^2) \times \mathrm{O}_2^{-}(q^2),\;
\mathrm{G}_2(q),\;
[2^9].\mathrm{SL}_3(2)\;(q=p),\;
P_{1,3,4}=[q^{11}].\mathrm{GL}_2(q)\,\mathrm{GL}_1(q)^2 .
\]

From \cite[Theorem 7 and Table 3]{AB2015} we know that \(H_0\) is large in \(G_0\) when it is of type \(\mathrm{G}_2(q)\), of type \(\mathrm{GL}_1^{\epsilon}(q) \times \mathrm{GL}_3^{\epsilon}(q)\) with \((q,\epsilon) \neq (3,+)\), of type \([2^9].\mathrm{SL}_3(2)\) with \(q=3\), or of type \(P_{1,3,4}\). In all these cases \(H_1\) is automatically large in \(G_1\).

Assume that $H_0$ is of type  $\mathrm{GL}_1^{\epsilon}(q)\times \mathrm{GL}_3^{\epsilon}(q)$ with $(q,\epsilon)= (3,+)$.
 By \cite[Table~8.50]{BHRD2013} we have \(|H_0| = |\mathrm{GL}_3(3)| =11232\) and \(|\mathcal{O}_1| = 24\). Then \(1< (|G_0| / |H_0|^3)^{1/2} < 2\), which implies that \(H_1\) is large in \(G_1\).
 Moreover, computation in {\sc Magma}~\cite{Magma} shows that $H_0=\mathrm{PSL}_3(3).2\cong \mathrm{Aut}(\mathrm{PSL}_3(3))$.

Assume that $H_0$ is of type $\mathrm{O}_2^{-}(q^2)\times \mathrm{O}_2^{-}(q^2)$.
 Using \cite[Table~8.50]{BHRD2013} we obtain $|H_0|=2^4(q^2+1)^2/d^2$, and  $|\mathcal{O}_{1} |=|\mathrm{Out}(G_0)|=24\log_p(q)/d^2$, where $d=(2,q-1)$,   
 Now $|H_0|\leq 2^4(q^2+1)^2\leq 2^4\cdot ( {5}/{4})^2q^4=25q^4$, and  $|\mathcal{O}_{1} |^2 \leq 24^2 \log_p(q))^2\leq 24^2\cdot 9q/8=648q$ (by Lemma~\ref{lm:q2f}).
By Lemma~\ref{lm:LUSO}, $q^{28}/8<|G_0|$.
From Eq.~\eqref{eq:H1} we derive  that $25^3\cdot 648\cdot  8>q^{23}$, which holds only for $q=2$.  For \(q=2\),  computation shows that \(1 < (|G_0| / |H_0|^3)^{1/2} < 2\), and so \(H_1\) is large in \(G_1\) although \(H_0\) itself is  not large in \(G_0\). Note that $H_0\cong \mathrm{D}_{10}^2.2.\mathrm{S}_2$ by \cite[Table~8.50]{BHRD2013}. 

Finally, let \(H_0\) be of type \([2^9].\mathrm{SL}_3(2)\) with \(q=p\) odd. Here \(H_0 \cong [2^9].\mathrm{SL}_3(2)\) and \(|\mathcal{O}_1| = 6\) by \cite[Table~8.50]{BHRD2013}. A computation shows that \(H_1\) is large in \(G_1\) precisely when \(q = 3\).  
\end{proof}


\subsection{Proof summary for Theorems~\ref{th:psl}--\ref{th:pso}}
We now consolidate the proofs of Theorems~\ref{th:psl}--\ref{th:pso}.
Recall that $G_0$ is a simple classical group, $H_0\in \mathcal{H}(G_0)$, $H_1=N_{\mathrm{Aut}(G_0)}(H_0)$ and $G_1=G_0H_1$. 
Lemma~\ref{lm:H1G1} asserts that $H_1$ is maximal in $G_1$, and $\mathcal{Q}(G_0,H_0)$ is nonempty if and only if $H_1$ is large in $G_1$.  
 
By Aschbacher's theorem (see e.g. \cite[Theorem 2.6.1]{BG2016}), every maximal subgroup of a classical group belongs to one of the subgroup collections \(\mathcal{C}_1, \dots, \mathcal{C}_8, \mathcal{S}, \mathcal{N}\). Our case-by-case analysis has covered all these possibilities:
\begin{itemize}
  \item For geometric subgroups (\(\mathcal{C}_1\)-\(\mathcal{C}_8\)), we treated linear groups in Subsection~\ref{sec:proof-psl}, unitary groups in Subsection~\ref{sec:proof-psu}, symplectic groups in Subsection~\ref{sec:proof-psp}, and orthogonal groups in Subsection~\ref{sec:proof-pso}.
  \item The  non-geometric subgroups in \(\mathcal{S}\) were handled in Subsection~\ref{sec:proof-as}.
  \item The small additional collection \(\mathcal{N}\) was analysed in Subsection~\ref{sec:proof-N}.
\end{itemize}
In each case we determined exactly when \(H_1\) is large in \(G_1\). Collecting the outcomes of all these subsections yields the complete lists of pairs \((G_0, H_0)\) presented in Theorems~\ref{th:psl}--\ref{th:pso}. \qed

\section*{Acknowledgments} 
The authors are grateful to the anonymous reviewers  for their valuable comments and suggestions.
  


\end{document}